\documentclass[a4paper,num-refs]{article}

\pdfoutput=1

\usepackage{IEK10} 
\usepackage{natbib}

\makeatother
\usepackage{amssymb}
\usepackage{xcolor}
\usepackage{mathtools, cuted}
\usepackage{tabularx}
\usepackage{eurosym}
\usepackage[utf8]{inputenc}

\newcommand{\vect}[1]{\mathbf{#1}}

\renewcommand{\min}[1][]{
	\ifthenelse{\isempty{#1}}{\operatorname{min}}{\ensuremath{\underset{#1}{\text{min}\,}}}
}

\newcommand*\circled[1]{\tikz[baseline=(char.base)]{
            \node[shape=circle,draw,inner sep=1pt] (char) {#1};}}
\newcommand{\RNum}[1]{\uppercase\expandafter{\romannumeral #1\relax}}

\def\IEK10{
  Forschungszentrum Jülich GmbH,
  Institute of Energy and Climate Research,
  Energy Systems Engineering (IEK-10),
  Jülich 52425,
  Germany
}
\def\RWTH{
  RWTH Aachen University
  Aachen 52062,
  Germany
}
\def\LTT{
  RWTH Aachen University,
  Institute of Technical Thermodynamics,
  Aachen 52056,
  Germany
}
\def\ETH{
  ETH Zürich,
  Energy \& Process Systems Engineering,
  Zürich 8092,
  Switzerland
}

\newcommand{\mytitle}{Simultaneous mixed-integer dynamic scheduling of processes and their energy systems}

\newcommand{\affil}{
  \begin{itemize}[leftmargin=3mm, itemsep=0mm]
    \item[$^a$]\IEK10
    \item[$^b$]\RWTH
    \item[$^c$]\ETH
    \item[$^d$]\LTT
  \end{itemize}
}

\def\firstAuthor{Florian Joseph Baader}
\newcommand{\myauthor}{\firstAuthor$^{a,b}$, André Bardow$^{a,c,d}$, Manuel Dahmen$^{a,*}$}
\newcommand{\correspondingauthor}{M. Dahmen, \IEK10 \newline E-mail: m.dahmen@fz-juelich.de}
\author{\myauthor}

\usepackage[
  colorlinks,
  linkcolor=blue,
  citecolor=blue,
  urlcolor=blue,
  pdftitle={\mytitle},
  pdfauthor={\firstAuthor}
]{hyperref}
\usepackage[capitalise, nameinlink]{cleveref}
\crefname{table}{Tab.}{Tab.}
\crefname{equation}{equation}{equations}

\newcommand{\setpgfexternalcounter}[1]{
  \makeatletter%
  \pgfkeysgetvalue{/tikz/external/figure name}\myexternalname
  \expandafter\gdef\csname c@tikzext@no@\myexternalname\endcsname{#1}%
  \makeatother
}

\begin{document}

  \thispagestyle{firststyle}

  \begin{center}
    \begin{large}
    {\fontsize{12}{14} \selectfont
      \textbf{\mytitle}}
    \end{large} \\
    \myauthor
  \end{center}

  \begin{footnotesize}
    \affil
  \end{footnotesize}

\begin{abstract}
Increasingly volatile electricity prices make simultaneous scheduling optimization desirable for production processes and their energy systems.
Simultaneous scheduling needs to account for both process dynamics and binary on/off-decisions in the energy system leading to challenging mixed-integer dynamic optimization problems.
We propose an efficient scheduling formulation consisting of three parts: a linear scale-bridging model for the closed-loop process output dynamics, a data-driven model for the process energy demand, and a mixed-integer linear model for the energy system. 
Process dynamics are discretized by collocation yielding a mixed-integer linear programming (MILP) formulation. 
We apply the scheduling method to three case studies: a multi-product reactor, a single-product reactor, and a single-product distillation column, demonstrating the applicability to multi-input multi-output processes.
For the first two case studies, we can compare our approach to nonlinear optimization and capture 82~\% and 95~\% of the improvement.
The MILP formulation achieves optimization runtimes sufficiently fast for real-time scheduling.
\end{abstract}


\newpage

\section*{Introduction}

Current efforts to reduce greenhouse gas emissions increase the share of renewable electricity production in many countries. 
Due to the intermittent nature of renewable electricity production, stronger volatility in electricity prices or even electricity availability is expected \citep{Merkert.2015}.
This price volatility may offer economic benefits to industrial processes that can dynamically adapt their operation and thus their power consumption in so-called demand response (DR)  \citep{Mitsos.2018}.
Ideally, demand response reacts to imbalances of electricity demand and supply and therefore also stabilizes the electricity grid \citep{Zhang.2016}.

A promising way to achieve DR is to consider volatile prices in scheduling optimization \citep{Merkert.2015} that determines the process operation for a time horizon in the order of one day \citep{Baldea.2014,Daoutidis.2018,Seborg.2010}.
However, industrial processes are often not supplied directly by the electricity grid but by a local on-site multi-energy system. 
The local multi-energy system supplies all energy demanded by the process, e.g., heating, cooling, or electricity, and exchanges electricity with the grid \citep{Voll.2013}.
Operating local energy systems is a complex task as these systems typically consist of multiple redundant units with non-linear efficiency curves and minimum part-load constraints leading to  discrete on/off-decisions \citep{Voll.2013}.
Thus, the electricity exchange between the energy system and the grid is not directly proportional to the process energy demand. Consequently, optimal DR scheduling must consider processes and their energy systems simultaneously.
Moreover, such a simultaneous scheduling can improve the efficiency of energy system operation by shifting process energy demand in time \citep{Bahl.2017}.
Still, scheduling is usually carried out sequentially: The process schedule is optimized first and only then the energy system operation is optimized \citep {MujtabaH.Agha.2010,Leenders.2019b}.

The simultaneous scheduling of processes and their energy systems leads to computationally challenging problems. Process scheduling can already be a very demanding task on its own if nonlinear process dynamics need to be considered \citep{FloresTlacuahuac.2010}; therefore, considering  dynamics is a major research topic in process systems engineering referred to as integration of scheduling and control \citep{Baldea.2014,Daoutidis.2018,Harjunkoski.2009,Engell.2012,Beal.2017,FloresTlacuahuac.2006}.
For DR problems, process dynamics are often scheduling-relevant \citep{Mitsos.2018,Baldea.2014,Daoutidis.2018,Caspari.2019,Otashu.2019} because the time to drive the process from one steady state to another steady state is often in the same order of magnitude as the electricity-price time steps.

\begin{figure}[h]
\centering
 \includegraphics{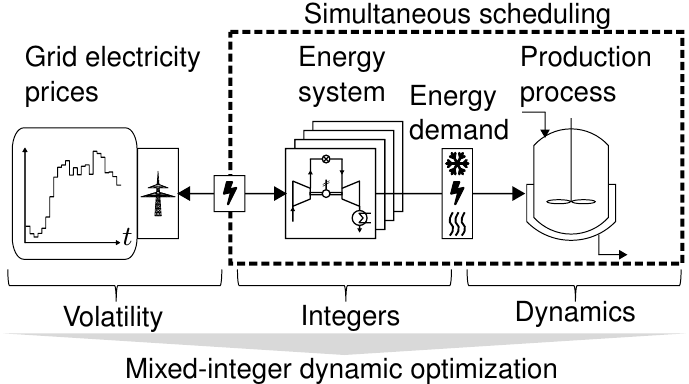}
  \caption{Volatile grid electricity prices call for a simultaneous scheduling of production processes and their local energy supply systems. While energy systems introduce integer decision variables, processes often exhibit scheduling-relevant dynamics. Simultaneous scheduling thus results in computationally challenging mixed-integer dynamic optimization (MIDO) problems.}
   \label{fig:problem}
\end{figure}

The desired simultaneous scheduling of processes and their energy systems is especially challenging due to the simultaneous presence of process dynamics and discrete on/off-decisions in the energy system (Figure \ref{fig:problem}). 
Because of the discrete decisions, standalone energy system optimization problems are preferably formulated as mixed-inter linear programs (MILPs) \citep{Voll.2013,MichaelJ.Risbeck.2017,Mitra.2013,Carrion.2006,SusanneSass.2020}.
A MILP formulation is usually applicable because: (i) Nonlinear part-load efficiencies can be approximated reasonably well using piece-wise affine functions \citep{SusanneSass.2020}, and (ii) the dynamics of the energy system units are negligible or can be captured using ramping constraints \citep{sass2019optimal}. 

As process dynamics are often scheduling-relevant, the simultaneous scheduling needs to be integrated with control. 
Even though conceptually, all approaches for the integration of scheduling and control can be used, the on/off-decisions significantly increase the computational complexity.
However, scheduling must be performed online.
\cite{IiroHarjunkoski.2014} state that generally optimization run times should be between 5 and 20 minutes.

In this work, we present a formulation for simultaneous scheduling of processes and their energy systems that aims at real-time-applicable runtimes.
We rely on two promising approaches from the integration of process scheduling and control: (i) dynamic scale-bridging models     \citep{Du.2015,Baldea.2015}, where the controlled process output is forced to follow a linear differential equation and (ii) dynamic data-driven models \citep{Mitsos.2018,Pattison.2016,Kelley.2018,Kelley.2018_2}.
Specifically, our formulation consists of three parts: \circled{1} a scale-bridging model considering the dynamics of the production process, \circled{2} a piece-wise affine dynamic data-driven model  for the energy demand of the process, and  \circled{3} a MILP energy system model with piece-wise affine approximations of nonlinear component efficiency curves. 
We discretize the linear differential equations in time using a high-order collocation scheme to receive linear constraints \citep{Biegler.2010}. 
Consequently, we achieve an MILP formulation for the entire scheduling problem.

A preliminary version of our approach has been presented in a conference contribution \citep{Baader.2020} where we considered DR for a building energy system. 
In the present contribution, we describe our method in more detail and apply it to three chemical production systems: a multi-product and a single-product continuous-stirred-tank-reactor (CSTR) both of which are cooled by three compression chillers, and a distillation column heated by two combined heat and power plants (CHPs) and an electrically-driven boiler.
The new method is explicitly compared against a standard sequential scheduling approach from industrial practice \citep{MujtabaH.Agha.2010}.
The remainder of this paper is structured as follows: In the second Section, the method is described in detail. In the third Section, a first case study considering a multi-product CSTR  is performed; in the fourth Section, a second case study considering a single-product CSTR is performed; and in the fifth Section, a third case study considering a distillation column is investigated.
The sixth Section concludes the work.

\section*{Method}
\label{sec:method}

\begin{figure}[h]
 \includegraphics[scale=1.0]{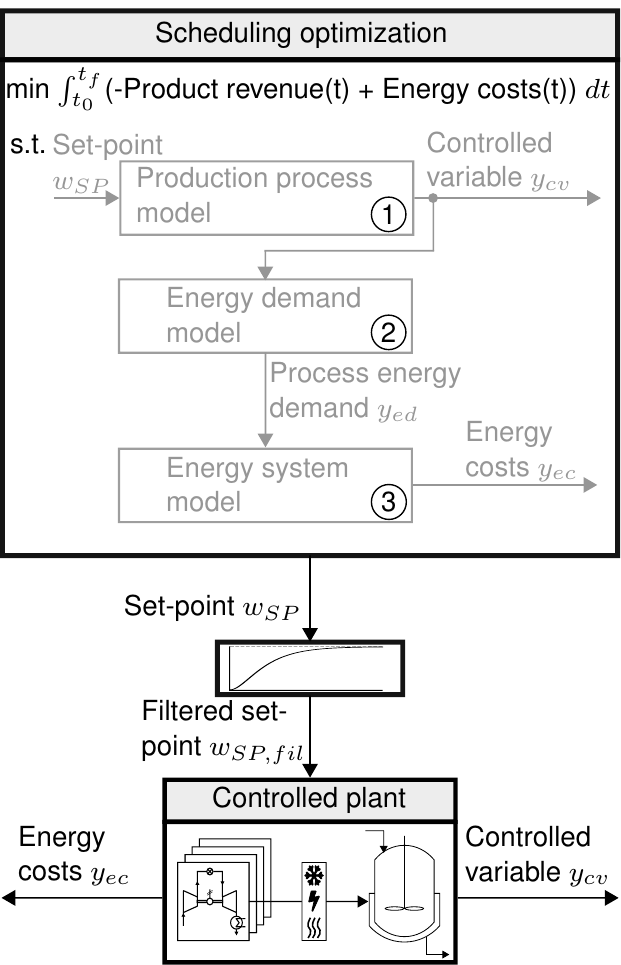}
 \centering
  \caption{Proposed simultaneous scheduling of processes and their energy systems based on our scheduling model consisting of three parts. A set-point filter converts the optimized piece-wise constant set-point $w_{SP}$  to a smooth filtered set-point $w_{SP,fil}$, which defines the desired linear closed-loop process behavior.}
   \label{fig:method}
\end{figure}

In this Section, we present our method for simultaneous dynamic scheduling of production processes and their energy systems. 
We refer to our method as \emph{simultaneous dynamic scheduling}.
The core of \emph{simultaneous dynamic scheduling} is an efficient scheduling model consisting of three parts: \circled{1} the production process,  \circled{2} the energy demand, and \circled{3} the energy system (Figure \ref{fig:method}). 
Model \circled{1} determines the controlled process output $y_{cv}$, e.g., the concentration in a reactor.
We use a scale-bridging model (SBM) proposed by Baldea and co-workers \citep{Du.2015,Baldea.2015} that describes a linear closed-loop response and represents the slow scheduling-relevant dynamics only. A linear SBM can be incorporated in scheduling optimization much more efficiently than a nonlinear full-order process model.
The SBM relies on an underlying control to enforce the desired linear closed-loop response.
The closed-loop response describes the evolution of the controlled variable $y_{cv}$ and its time derivatives depending on the set-point $w_{SP}$:
\begin{align}
    \label{eq:y_cv_equation}
    y_{cv} + \sum_{i=1}^{r}\tau_i \frac{d^i y_{cv}}{dt^i}=w_{SP}
\end{align}
In \cref{eq:y_cv_equation}, $r$ is the order of the SBM and $\tau_i$ are time constants. We discuss both order and time constants in the following subsection. 
To linearize the closed-loop response, we propose to place a set-point filter \citep{Corriou.2018} in front of the controlled plant (Figure \ref{fig:method}). 
This set-point filter converts the piece-wise constant set-point $w_{SP}$ given by the scheduling optimization to a smooth filtered set-point $w_{SP,fil}$ that can be tracked by the underlying process control such that $y_{cv} \approx w_{SP,fil}$.
In essence, we assume that the linear dynamics of the set-point filter can model the process output dynamics for the scheduling-relevant time-scale.
Instead of the combination of set-point filter and tracking control, previous publications used exact input-output feedback linearization \citep{Du.2015, Daoutidis.1992, Corriou.2018} or scheduling-oriented model predictive control (SO-MPC) \citep{Baldea.2015}. 
The proposed set-point filter increases the flexibility of the scale-bridging approach as it allows to use non-model-based tracking controls, e.g., PID-control \citep{Corriou.2018}, as well. 

A disadvantage of the SBM is the resulting conservatism because the time constants need to be chosen such that the desired linear closed-loop response can be realized in all operating regimes. 
Therefore, the closed-loop response is slower than necessary in some operating regimes.

The main advantage is that the scale-bridging \cref{eq:y_cv_equation} is more than an approximation: Whenever the actual value of $y_{cv}$ deviates from the closed-loop response described by \cref{eq:y_cv_equation}, the underlying control acts to bring the controlled variable $y_{cv}$ back to the desired closed-loop trajectory.
Consequently, deviations of the controlled variable  from its optimized trajectories are kept small.

Note that, in this paper, we study the case of a single scale-bridging model, whose application is straight-forward for single-input single-output (SISO) processes where the scale-bridging model deals with the only controlled variable $y_{cv}$. 
In the multi-input multi-output (MIMO) case, a vector of output variables $\vect{y}_{cv}$ is controlled; still, the number of slow scheduling-relevant variables is typically small \citep{Baldea.2014,Baldea.2015,Pattison.2016}.
Moreover, for flexible DR operation, often, only one scheduling-relevant quantity $\rho$ is varied to shift energy demand in time, e.g., the production rate or product purity. 
A single scale-bridging model for this quantity $\rho$ is sufficient if all controlled outputs are either maintained constant irrespective of the scheduling-relevant quantity $\rho$ or are alternatively coupled with $\rho$.
For instance, the hold-up of process units might be maintained constant irrespective of flexible operation.
An example for coupling of controlled process outputs is given in the initial scale-bridging model paper where  \cite{Du.2015} consider a reactor with concentration and temperature as controlled outputs but only a scale-bridging model of the concentration during scheduling.
After scheduling, \cite{Du.2015} derive a set-point signal for the temperature from the set-point signal of the concentration. 
Consequently, often, only one scale-bridging model needs to be tuned even though multiple inputs and outputs are present.
Such a case is also shown in our third case study, where we consider a 4$\times$4 MIMO process.

Model \circled{2} is a dynamic data-driven model \citep{Mitsos.2018}  that determines the process energy demand $y_{ed}$ taking the current state of the production process as inputs, i.e., the controlled variable $y_{cv}$ and its time derivatives.
In principle, a wide range of dynamic data-driven models derived from recorded data, or mechanistic models can be used here. 
Examples of dynamic data-driven models being applied successfully in dynamic demand response optimization can be found in the literature \citep{Pattison.2016,Kelley.2018,Kelley.2018_2,Tsay.2019}. 
Our energy demand model \circled{2} can be dynamic and mixed-integer but must be linear as we aim for an MILP formulation. 
Note that the energy demand can, in general, not only depend on the controlled outputs but also on other uncontrolled states. 
In such cases, the data-driven energy demand model needs to have internal states that approximate the internal dynamics of the real process. 
Models with internal states, such as Hammerstein-Wiener models, are common in system identification and have also been used in  demand-response applications \citep{Kelley.2018_2}. An energy demand model with an  internal state is demonstrated in our third case study.

In contrast to the controlled process output $y_{cv}$, deviations between the actual process energy demand $y_{ed}$ and the model prediction are not corrected.
Instead, we assume that such deviations are compensated by the energy system, which is reasonable if (i) the energy system can react significantly faster than the process and (ii) the energy system has spare capacity larger than the maximum error of the data-driven model.

In this paper, we assume that the energy demand is the only uncontrolled process quantity that is relevant for the scheduling objective function. 
In principle, other  uncontrolled quantities could be relevant for the scheduling objective function as well.
For instance, raw material consumption could vary due to flexible operation and thus cause additional costs that should not be neglected in the optimization.
In such cases, data-driven models for those quantities need to be derived and added in the same way as for the energy demand.

Model \circled{3} is the energy system model that determines the energy costs depending on the energy demand. 
The structure of the energy system is modeled by energy balances that connect the energy system components with demands.
Moreover, the efficiency of individual energy system components is modeled as a function of
the part-load fraction.
Thus, the required input power $P_{c,in}$ of an energy system component $c$ is a nonlinear function of the desired output power $P_{c,out}$ \citep{SusanneSass.2020}.
To obtain an MILP formulation, we follow the established approach of modeling part-load efficiency curves as piece-wise affine functions \citep{Voll.2013}.
In general, piece-wise affine efficiency curves require binary variables.
Binary variables increase the computational burden; however, they can be avoided if the input power is a convex function of the output power \citep{Carrion.2006}, which is the case for many energy system components of practical relevance \citep{Voll.2013}.

By combining the three models \circled{1} - \circled{3}, we receive a linear differential algebraic equation system (DAE) containing integers:
\begin{align}
    \label{eq:gen_SBM_1}
     \frac{d\vect{x}}{dt} &= \vect{f}(\vect{x},\vect{y},\vect{z},\vect{w}_{SP})  = \vect{A}\vect{x} + \vect{B}\vect{y} + \vect{C}\vect{z} + \vect{D}\vect{w}_{SP} \\
     \label{eq:gen_SBM_2}
      \vect{0} &= \vect{g}(\vect{x},\vect{y},\vect{z},\vect{w}_{SP})  = \vect{E}\vect{x} + \vect{F}\vect{y} + \vect{G}\vect{z} + \vect{H}\vect{w}_{SP} 
\end{align}
In \cref{eq:gen_SBM_1,eq:gen_SBM_2}, $\vect{x}$ are the differential states, $\vect{y}$ are continuous variables, $\vect{z}$ are discrete variables, $\vect{w}_{SP}$ are set-points, $t$ is time, $\vect{f}$ and  $\vect{g}$ are functions that are linear in $\vect{x},\vect{y},\vect{z},\vect{w}_{SP}$, and  $\vect{A}-\vect{H}$ are matrices.
Note that all variables are functions of time although not stated explicitly to improve readability.

We choose a discrete-time MILP formulation for our simultaneous scheduling problem because 
in case of variable electricity prices, discrete-time formulations usually perform better than continuous-time formulations \citep{Castro.2009}, as the electricity markets imposes a discrete time structure, e.g., hourly constant prices.
As our model consists of linear differential equations,  time discretization with collocation in discrete time leads to linear constraints \citep{Biegler.2010}.

In the following, we discuss the scale-bridging model parameters and the scheduling optimization problem.

\subsection*{Scale-bridging model parameters}
\label{sec:met_filter}

For the scale-bridging model \circled{1}, we have to determine the order $r$ and the time constants $\tau_i$ from \cref{eq:y_cv_equation}, as well as upper and lower bounds for the set-point $w_{SP}^{max}$, $w_{SP}^{min}$, respectively.
The order $r$ should be chosen such that the resulting reference trajectory for the controlled variable $y_{cv}$ can be realized by the process. 
For instance, if a process reacts with second-order dynamics to input changes, a first-order set-point filter is not reasonable.
Thus, the order $r$ should reflect how many stages of inertia a change in the manipulated variable $u$ has to overcome before changing the controlled variable $y_{cv}$. 
If a process model is available, the order $r$ can be derived mathematically by analyzing the relative degree of the process model defined as the number of times the controlled variable $y_{cv}$  
has to be differentiated with respect to time until the manipulated variable $u$ appears explicitly \citep{Corriou.2018}.
If no process model is available, the order $r$ needs to be chosen based on knowledge or intuition about the main inertia of the process.
In this way, the set-point filter reflects the characteristics of the open-loop system.
The employed controller might add additional dynamics to the closed-loop system, for instance the integral part of a PI-controller.  
Using the relative order of the open-loop process to choose the relative order of the set-point filter is thus only reasonable if it can be assumed that either the controller dynamics are significantly faster than the process dynamics or the proportional part of the controller dominates over the integral part.
That is, our approach might not work well in a case where the controller adds significant dynamics to the closed-loop system. 
In such a case, the controller can be changed or alternative approaches that explicitly model the control behavior might be used instead. 
For instance, \cite{LisiaS.Dias.2018} optimize a schedule using a sequential optimization approach with the controller being part of the simulation model. 
Thus, they optimize the closed-loop system.
\cite{JeromeE.J.Remigio.2020} explicitly account for the behavior of an underlying model predictive controller (MPC) by adding the KKT-conditions of the MPC to the scheduling optimization problem.

As discussed by \cite{Baldea.2015}, the choice of the time constants $\tau_i$ is critical for the performance of the scale-bridging approach: 
On the one hand, if the time constants are too small, the scale-bridging dynamics are too fast and cannot be realized by the controlled process.
On the other hand, if the time constants are too large, the scale-bridging dynamics are overly conservative and process flexibility is wasted.
However, a rigorous way to tune the time constants $\tau_i$ is missing in the literature.

In this paper, we argue that the time constants $\tau_i$ need to be tuned simultaneously with the set-point bounds $w_{SP}^{max}$ and $w_{SP}^{min}$.
For illustration, we consider a transition of the controlled variable starting from a small value $y_{cv}^{start}$ and ending at a new steady state with a higher value $y_{cv}^{end}$ close to the maximum allowable value $y_{cv}^{max}$.
To 
speed up the transition,
scheduling optimization might choose a set-point $w_{SP,elevated}$ which is elevated above $y_{cv}^{end}$ and even $y_{cv}^{max}$ for a certain period of time. 
However, choosing an elevated set-point value can lead to dynamics that are too fast for the controlled process.
In particular, if the time constants $\tau_i$ are small, the scale-bridging dynamics are already fast and an elevated set-point may drive the controlled variable to infeasible values.
A trade-off arises because we want to choose small time constants in general but also want to avoid slow transitions towards the bounds of the controlled variable.

In our case studies, we tune the scale-bridging parameters using a simple heuristic relying on simulations. 
Alternatively, existing knowledge about the time constants of the process, or measurements can be used to calibrate the SBM.

\subsection*{Scheduling optimization problem}
To derive a complete problem formulation based on \cref{eq:gen_SBM_1,eq:gen_SBM_2}, we add a suitable objective function, discretize time, and add inequality constraints to account for variable bounds, minimum part-load, and problem specific constraints.

The objective $\Phi$ in simultaneous DR scheduling is to maximize cumulative product revenue $\Phi_{Product}$ at final time $t_f$ minus the cumulative energy costs $\Phi_{Energy}$ at final time:
\begin{align}
    \label{eq:obj}
    &\text{min}~\Phi = - \Phi_{Product}(t_{f}) + \Phi_{Energy}(t_{f}) \\
    \label{eq:obj_product}
    &\frac{ d\Phi_{Product}}{dt} = \sum_{p \in \mathbb{P}} K_p q_p\\
    \label{eq:obj_energy}
    &\frac{ d\Phi_{Energy}}{dt} = \sum_{e \in \mathbb{E}} K_{e} P_{e}\\ \nonumber
    &\text{with}~~\Phi_{Product}(t_0)=\Phi_{Energy}(t_0)=0
\end{align}
Here, $\mathbb{P}$ is the set of products, $q_p$ the flow rate of product $p$, and $K_p$ the price of $p$. Similarly, $\mathbb{E}$ is the set of end-energy forms consumed, $K_{e}$ is the time-dependent price of energy $e$, and $P_{e}$ is the consumed power of energy $e$.  $t_0$ denotes the initial time.

\begin{figure}[h]
\centering
    \includegraphics{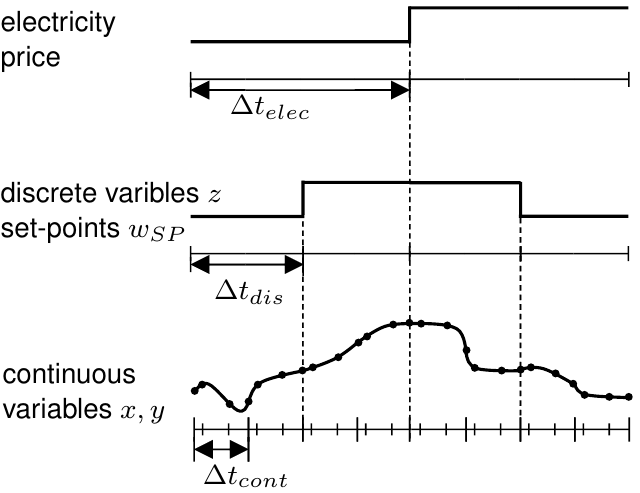}
    \caption{Three time grids used for discretization with timesteps $\Delta t_{elec}$, $\Delta t_{dis}$, $\Delta t_{cont}$, respectively. }
    \label{fig:time_grids}
\end{figure}

For time discretization, we use three time grids  (Figure \ref{fig:time_grids}).
Grid 1 is given by the electricity market and contains piece-wise constant electricity prices with time step $\Delta t_{elec}$, e.g., 1~h or 15~min. 
Grid 2 contains discrete decision variables $\vect{z}$ and piece-wise constant set-points $\vect{w}_{SP}$. 
The resolution of grid 2 should not be too fine as it increases the number of integer variables and thus the combinatorial complexity of the optimization problem.
Still, it should be possible to alter discrete decisions $\vect{z}$ and set-points $\vect{w}_{SP}$ at least at every step change of electricity prices.
Thus, the electricity price time step resolution should constitute a lower bound on the resolution of grid 2.
Making grid 2 finer than grid 1 by selecting time steps $\Delta t_{dis} < \Delta t_{elec}$, gives a higher flexibility and thus might enable higher profits.
We recommend to use time steps with lengths $\Delta t_{dis} = \frac{1}{n_1}\Delta t_{elec}$ with $n_1$ being a small natural number.

Grid 3 is used for continuous variables $\vect{x}, \vect{y}$.
Differential states $\vect{x}$ are discretized using collocation. 
Similar to the argument above, we propose to use finite elements with length $\Delta t_{cont} = \frac{1}{n_2}\Delta t_{dis}$. 
The natural number $n_2$ is chosen to be greater than or equal to one because whenever electricity prices, discrete variables, or set-points perform step-changes, a new collocation element is necessary such that non-smoothness in differential states $\vect{x}$ is possible. 
In result, $\vect{x}$ are continuous at the border of collocation elements but first derivatives are allowed to perform step changes.

Within a finite element $f_e$ of grid 3, a collocation polynomial $x_{f_e}$ of order $N_{cp}$ is used to discretize differential states \citep{Biegler.2010,Nicholson.2018}: 
\begin{align}
    \label{eq:col_1}
    x_{f_e} (\overline{\tau}) &= \sum_{j=0}^{N_{cp}}l_j(\overline{\tau})x_{f_e,j},  ~~~~~    \overline{\tau} \in [0,1]\\
    \label{eq:col_2}
    l_j&=\prod_{k=0,k\neq j}^{N_{cp}} \frac{\overline{\tau} - \overline{\tau}_k}{\overline{\tau}_j-\overline{\tau}_k}\\
    \label{eq:x_dot_col}
    \frac{dx}{dt}\bigg|_{t_{f_e,k}} &= \frac{1}{\Delta t_{cont}} \sum_{j=0}^{N_{cp}}x_{f_e,j}\frac{dl_j(\overline{\tau}_k)}{d\overline{\tau}}
\end{align}
In \cref{eq:col_1,eq:col_2}, the $l_j$ are Lagrange basis polynomials, $\overline{\tau}$ is the scaled time within a finite element, and $x_{f_e,j}$ are state values at discretization points. In \cref{eq:x_dot_col}, $\frac{dx}{dt}\bigg|_{t_{f_e,k}}$ is the approximated time derivative at a collocation point $k$, which is set equal to the right hand side of the linear differential \cref{eq:gen_SBM_1} for every time point ${t_{f_e,k}}$. The term $\frac{dl_j(\overline{\tau}_k)}{d\overline{\tau}}$ is a constant parameter in the optimization because it only depends on $\overline{\tau}$. Moreover, as we choose discrete time, $\Delta t_{cont}$ is constant.
Therefore, $x_{f_e,j}$ are the only optimization variables, and thus, discretization with \cref{eq:x_dot_col} leads to linear constraints.

As inequality constraints, we consider upper and lower bounds for all variables, minimum part-load constraints for energy system components, and problem-specific constraints, e.g., minimum production targets. Minimum part-load constraints  are realized with a binary variable $z_{on,c}$ that indicates if the output power $P_{c,out}$ of an energy system component $c$ is zero or between the minimal and maximal allowed value, $P_{c,out}^{min}$ and $P_{c,out}^{max}$, respectively: 
\begin{align}
    \label{eq:min_part_load}
    z_{c,on}P_{c,out}^{min} \leq P_{c,out} \leq z_{c,on}P_{c,out}^{max}
\end{align}

Assembling the discussed equations gives the simultaneous scheduling optimization problem for the production process and its energy system.

Our discussion focuses on energy-intensive processes that can be operated flexibly and exhibit scheduling-relevant dynamics.
The energy demand of other processes also present at the same chemical production site can be integrated in our MILP scheduling optimization problem in straight-forward manner: 
The energy demand of processes with no or negligible flexibility can be expressed as predefined time-varying demands that can be added to the energy balances as constant terms.
For processes with negligible dynamics, quasi-steady state can be assumed and the steady-state energy demand can be modeled as a piece-wise affine function \citep{Schafer.2020}.

\section*{Case study 1: Multi-product reactor}
In this Section, we assess the computational performance of our method in a first  case study considering a multi-product reactor. 
We benchmark the economic value of simultaneous dynamic scheduling to a standard sequential scheduling and to a nonlinear scheduling optimization with the true process model.
\subsection*{Setup}

\begin{figure}[h]
\centering
    \includegraphics{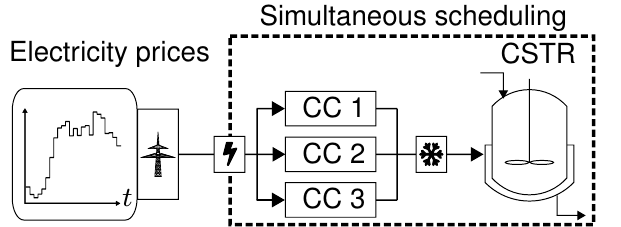}
    \caption{Case study 1: Simultaneous scheduling of a continuous-stirred-tank-reactor (CSTR) cooled by three compression chillers (CC1, CC2, CC3). Time-varying electricity prices provide an economic incentive for DR. }
    \label{fig:case_study}
\end{figure}

The setup of the case study is visualized in Figure \ref{fig:case_study}.
An exothermic multi-product CSTR is cooled with three compression chillers (CC). 
We use an exemplary reactor model from  \cite{Petersen.2017}. 
In the CSTR, a component A reacts to a component B. 
The reactor can produce three products $\RNum{1}$, $\RNum{2}$, $\RNum{3}$, which are defined by the desired concentration of component A, $C_A$.
We assume a small tolerance of +/- 0.01~$\frac{mol}{L}$ such that for each product, we obtain a product band. 
Whenever the concentration $C_A$ is within one of the three product bands, the associated product is produced. 
If the concentration is outside of the three product bands, which happens necessarily during transitions, no product is produced. 
For illustration, we consider prices of 1, 0.75, and 0.5 money unit (MU) (Table \ref{tab:products}) and require a minimum daily production of 5 hours and a maximum daily production of 8 hours for each product.

In the CSTR model, the rate of change for the concentration of component A is given by the material balance and the rate of change for the temperature $T$ is given by the energy balance:

\begin{align}
    \label{eq:case_study_1}
    &\frac{dC_A}{dt}=\frac{q}{V}(C_{A,feed}-C_A)-ke^{-\frac{E_A}{RT}}C_A\\
    \label{eq:case_study_2}
    &\frac{dT}{dt}=\frac{q}{V}(T_{feed}-T)-\frac{k \Delta H_r}{\varrho C_P}e^{-\frac{E_A}{RT}}C_A - \frac{Q_{cool}}{\varrho C_P V}
\end{align}
In \cref{eq:case_study_1,eq:case_study_2}, $q$ is the flow rate, $V$ the reactor volume, $C_{A,feed}$  the feed concentration, $k$ the reaction constant, $E_A$ the activation energy, $R$ the gas constant, $T_{feed}$ the feed temperature, $\Delta H_r$ the enthalpy of reaction, $\varrho$ the density, $C_P$ the heat capacity, and $Q_{cool}$ is the cooling provided to the reactor.
The parameter values listed in Table~S1 of the Supplementary information (SI) are exemplary values from \cite{Petersen.2017}, except that we varied the activation energy $E_A$ to obtain an operating temperature where cooling with compression chillers is a realistic option.

An efficient chiller is used for base-load cooling, whereas  chiller 2 has a medium coefficient of performance (COP), and chiller 3, which has a low COP, is used for peak cooling (see Table \ref{chiller}).
We use the compression chiller model from \cite{Voll.2013} with a minimum part load of 20~\% and a coefficient of performance depending on nominal COP, $COP_{CC,i}^{nom}$, cooling load $Q_{CC,i}$, and nominal cooling load $Q_{CC,i}^{max}$:
\begin{align}
    \label{eq:case_study_3}
    &COP_{CC,i} = COP_{CC,i}^{nom}( 0.8615 q_{CC,i}^3 - 3.5494 q_{CC,i}^2     + 3.679q_{CC,i}+0.0126  ),    \text{with}~ q_{CC,i}=\frac{Q_{CC,i}}{Q_{CC,i}^{max}}
\end{align}

Note that this case study is meant to be an illustrative example rather than a real case.
It allows us to study whether the proposed method is able to consider process dynamics and discrete on/off-decisions for energy system components simultaneously. 
We want to stress that even though the original nonlinear process model is a small-scale model, the resulting scale-bridging model (\circled{1}) would have the same basic structure and computational complexity if a larger process model would be considered as the number of scheduling-relevant dynamics is typically small \citep{Baldea.2014,Baldea.2015,Pattison.2016}.

\begin{table}[bt]
\centering
    \caption{Product band $\left[C_{A,p}^{min}, C_{A,p}^{max}\right]$ in mol/L, price $K_p^P$ in money unit (MU)/$\text{m}^3$, and cooling power in steady state production
     $Q_{cool,p}^{steady}$ in $\frac{MJ}{h}$, 
    for products $p \in \{\RNum{1} , \RNum{2}, \RNum{3}\}$}
    \large
    \renewcommand\theadfont{\large}
    \label{tab:products}
    \begin{tabular}{cccc}
    \hline
    \thead{$p$} & \thead{$\left[C_{A,p}^{min},C_{A,p}^{max}\right]$} & \thead{$K_p^P$   } & \thead{$Q_{cool,p}^{steady}$  }\\
    \hline
    $\RNum{1}$ &  [0.09, 0.11] & 1 & 6.05\\
    $\RNum{2}$ & [0.29, 0.31] & 0.75 & 5.43\\
    $\RNum{3}$ & [0.49, 0.51] & 0.5 & 4.65\\
    \hline  
    \end{tabular}
\end{table}

\begin{table}[bt]
\centering
    \caption{Nominal cooling power $Q_{CC,i}^{max}$ and coefficient of performance $COP_{CC,i}^{nom}$ for compression chillers}
    \large
    \renewcommand\theadfont{\large}
    \label{chiller}
    \begin{tabular}{ccc}
    \hline
    \thead{compression chiller} & \thead{$Q_{CC,i}^{max}$ [$\frac{MJ}{h}$]} & \thead{$COP_{CC,i}^{nom} [-]$} \\
    \hline
    1 & 4.8 & 6\\
    2 & 2.3 & 4.5\\
    3 & 1.5 & 3\\
    \hline  
    \end{tabular}
\end{table}

We employ conventional PID control \citep{Corriou.2018} to track the filtered set-point for the concentration $C_A$ by manipulating the cooling power $Q_{cool}$:
\begin{align}
    \label{eq:PID}
    Q_{cool} = K_P *\left(e +\tau_D \frac{de}{dt} + \frac{1}{\tau_I} \int_0^t e dt\right) +Q_{cool}^{PID,0}~~~\text{, with}~ e=w_{SP,fil} - C_A
\end{align}
The controller parameters in \cref{eq:PID} are: $K_P$, $\tau_D$, $\tau_I$, and $Q_{cool}^{PID,0}$. 
We choose $Q_{cool}^{PID,0}$ to be the steady-state cooling power of product $\RNum{2}$ (Table \ref{tab:products}) and manually tune the remaining controller parameters in a simulation such that the filtered set-point $w_{SP,fil}$ is tracked stable and accurately. The resulting parameters are: $K_P = 1000~\frac{MJ~L}{h~mol}$, $\tau_D=0.1~h$, and $\tau_I=0.2~h$. 
The stable and accurate set-point tracking is shown in the following (Figure \ref{fig:transitions}). 

\subsection*{Simultaneous dynamic scheduling}
To apply our simultaneous dynamic scheduling method, we now set up the three parts of our model and the scheduling optimization problem as presented in the Method Section. 

\subsubsection*{Scale-bridging production process model}
\label{sec:filter_tuning}
\begin{figure}
\centering
    \includegraphics{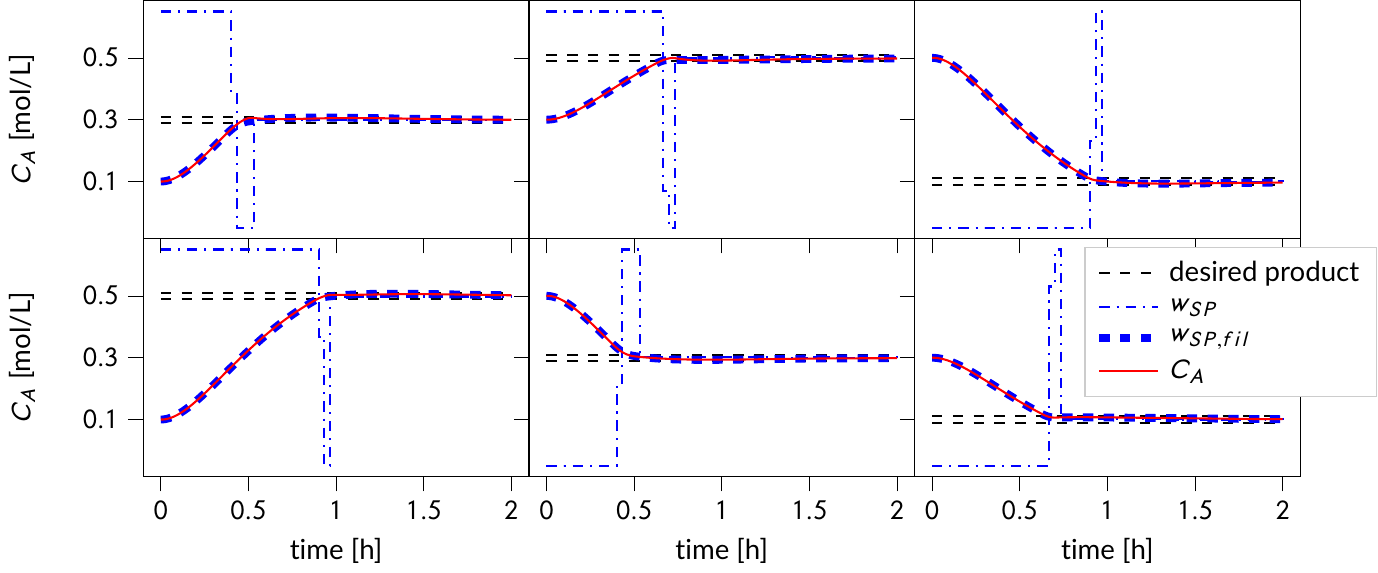}
    \caption{Six possible transitions between the three products with the chosen time constant $\beta=0.36~h$ and set-point elevation  $w_{SP}^{elevation}=0.15~\frac{mol}{L}$ (compare to Figure 5). The piece-wise constant set-point from optimization $w_{SP}$ results in a filtered set-point $w_{SP,fil}$, which can be tracked accurately such that the actual value of the concentration $C_A$ resulting from the controlled nonlinear process model matches the filtered set-point $w_{SP,fil}$ well.}
   \label{fig:transitions}
\end{figure}

As discussed in the Method Section,
we need to choose the order $r$, the time constants $\tau_i$, and the set-point bounds $w_{SP}^{min}$, $w_{SP}^{max}$ for the scale-bridging production process model \circled{1}.
The relative order of the open-loop process is 2, as can be seen from the physical process model: The manipulated variable $Q_{cool}$ does not appear in the first derivative of the controlled variable $C_A$ (\cref{eq:case_study_1}). 
If \cref{eq:case_study_1} is differentiated with respect to time and the term $\frac{dT}{dt}$ is replaced using \cref{eq:case_study_2}, the second time derivative $\frac{d^2 C_A}{dt^2}$ appears as an explicit function of the input $Q_{cool}$. Thus, the open-loop system has the relative degree 2.
The more descriptive explanation is that a change in the cooling power $Q_{cool}$ has to first overcome the inertia of the temperature $T$ and then the inertia of the concentration $C_A$.
As discussed in the Method Section
, we choose the order of the set-point filter equal to the relative order of the open-loop process, i.e., $r=2$. 
Our simulation results show that the second-order response can in fact be realized by the closed-loop system (Figure \ref{fig:transitions}).

Second-order systems are described in control theory by the time constant of their natural oscillation $\beta$ and a damping coefficient $\zeta$ \citep{Corriou.2018}.
The two tunable time constants $\tau_1$ and $\tau_2$ can be expressed as:
    \begin{align}
        \tau_1 = 2\zeta\beta \\
        \tau_2 = \beta^2
    \end{align}
Following \cite{Du.2015}, we choose a critically damped response, i.e., $\zeta = 1$, as we want to have fast but no oscillating dynamics.

In the following, we describe the heuristic procedure used to define the remaining time constant $\beta$ simultaneously with the bounds for the set-point $w_{SP}^{max}$ and $w_{SP}^{min}$.
The allowed range of the set-point must at least cover the operating range of the concentration $C_A$ which is between $C_A^{min} = 0.1~\frac{mol}{L}$ and $C_A^{max} = 0.5~\frac{mol}{L}$.
However, 
as discussed in the Method Section, it is reasonable to allow elevated set-points in order to avoid overly conservative transitions towards the bounds of the concentration $C_A$.
We introduce an elevation constant $w_{SP}^{elevation}$ and calculate the bounds of the set-point as:
    \begin{align}
        w_{SP}^{max} = C_A^{max} + w_{SP}^{elevation} \\
        w_{SP}^{min} = C_A^{min} - w_{SP}^{elevation}
    \end{align}
We want to find a combination of $\beta$ and $w_{SP}^{elevation}$ that (i) is feasible, i.e., the filtered set-point can be tracked accurately without oscillations, and (ii) allows for fast product transitions. In the following, we first present a routine to evaluate the feasibility and speed for a given combination of $\beta$ and $w_{SP}^{elevation}$ and then explain how we explore the space of possible combinations.

To evaluate a combination of $\beta$ and $w_{SP}^{elevation}$, we first optimize and then simulate all 6 possible transitions between the 3 products bands.
Thus, each of the 6 transitions starts at a product $p_i$ and ends at another product $p_j$, $i\neq j$.
For a given combination of $\beta$ and $w_{SP}^{elevation}$, we perform the following four steps for each transition:
    \begin{enumerate}
        \item As we assume that a fast transition is more critical than a slow one, we optimize a trajectory of set-points $w_{SP}(t)$ using model \circled{1} to start from product $p_i$ and reach the product band of product $p_j$ as fast as possible. 
       To generate this as-fast-as-possible set-point trajectory, we use exactly the same constraints as in the scheduling optimization. The resulting optimization problem reads:
        \begin{align}
            \label{eq:obj_filter_tuning}
            &\hspace{3.0cm}  
            \underset{w_{SP}(t)}{\text{min}}~~ - \sum_{t_0}^{t_f} z_{p_j}(t) \\
            \text{s.t. }~~& w_{SP,fil}(t) + 2\beta\frac{d w_{SP,fil}}{dt}\bigg|_t + \beta^2\frac{d^2 w_{SP,fil}}{dt^2}\bigg|_t = w_{SP}(t)~~~&\forall t \in \left[t_0, t_f\right] \label{eq:tuning_SBM}\\
            \label{eq:z_p_1}
            &w_{SP,fil}(t) - (C_{A,p_j}^{min} + \varepsilon_p) \geq -(1-z_{p_j}(t))  ~~~~~ &\forall t \in \left[t_0, t_f\right] \\
            \label{eq:z_p_2}
            &(C_{A,p_j}^{max} - \varepsilon_p)-w_{SP,fil}(t)  \geq -(1-z_{p_j}(t))  ~~~~~ &\forall t \in \left[t_0, t_f\right] \\
            \label{eq:tuning_set_point_bounds}
            &C_A^{min} - w_{SP}^{elevation} \leq w_{SP}(t) \leq C_A^{max} + w_{SP}^{elevation} &\forall t \in \left[t_0, t_f\right] \\
            \label{eq:initial_wsp}
            &w_{SP,fil}(t_0) = C_{A,p_i}\\
            \label{eq:initial_dwsp_dt}
            &\frac{d w_{SP,fil}}{dt} \bigg|_{t_0} = 0
        \end{align}
        In \cref{eq:obj_filter_tuning,eq:tuning_SBM,eq:z_p_1,eq:z_p_2,eq:tuning_set_point_bounds,eq:initial_wsp,eq:initial_dwsp_dt}, $z_{p_j}(t)$ is a binary indicating if the filtered set-point $w_{SP,fil}$ has reached the band of the desired product $p_j$. Accordingly, $z_{p_j}(t)$ is 1 if $C_{A,p_j}^{min} + \varepsilon_p \leq w_{SP,fil}(t) \leq C_{A,p_j}^{max} - \varepsilon_p$ where $\varepsilon_p$ is a small tolerance which we set equal to 0.003. The value of $z_{p_j}(t)$ is enforced by \cref{eq:z_p_1,eq:z_p_2}.  \Cref{eq:tuning_SBM} is the SBM (\circled{1}) and \cref{eq:tuning_set_point_bounds} bounds the piece-wise constant set-point $w_{SP}$. Time discretization with collocation converts the dynamic optimization problem to an MILP.
        \item We take the resulting set-point trajectory $w_{SP}(t)$ as input to a simulation of the set-point filter, the underlying PID-control, and the nonlinear process model.
        \item Based on the simulation result, we check if the transition is feasible. We consider a transition to be feasible if (i) the concentration reaches the product band of the desired product $p_j$ and (ii) the concentration stays inside the band  of $p_j$ once this product band is reached.
        \item We store the time needed to reach the product band $\Delta t_s$, which is calculated from the simulation result as the simulation gives the true closed-loop response. 
    \end{enumerate}
A combination of $\beta$ and $w_{SP}^{elevation}$ is considered feasible if all 6 transitions are feasible.
We measure the quality of feasible parameter combinations by the sum of all 6 transition times, i.e.,
    \begin{align}
        \Delta t_{sum} = \sum_{s \in \mathbb{S}} \Delta t_s,
    \end{align}
    where $\mathbb{S}$ is the set of possible transitions.
As we aim for fast transitions, we prefer feasible combinations of $\beta$ and $w_{SP}^{elevation}$ with a small value of $\Delta t_{sum}$.

To explore the space of possible combinations, we first set the set-point elevation to zero, i.e., $w_{SP}^{elevation} = 0~\frac{mol}{L}$, and start with $\beta = 0.1~h$.
We increase $\beta$ by steps of size $\Delta \beta = 0.01~h$ and evaluate if all six product transitions are feasible.
That is, for every value of $\beta$, we repeat steps 1 - 4 for all six transitions.
The smallest  time constant $\beta$ for which all 6 transitions are feasible is $\beta^{min} = 0.26~h$.
Starting from $\beta^{min} = 0.26~h$, we continue to increase $\beta$ and additionally allow a set-point elevation.
For every value of $\beta$, the set-point elevation $w_{SP}^{elevation}$ is increased by steps of size $\Delta w_{SP}^{elevation} = 0.01~\frac{mol}{L}$ until one of the transitions becomes infeasible. 
That way, we find the highest possible set-point elevation for every $\beta$.

\begin{figure}   
\centering
   \includegraphics{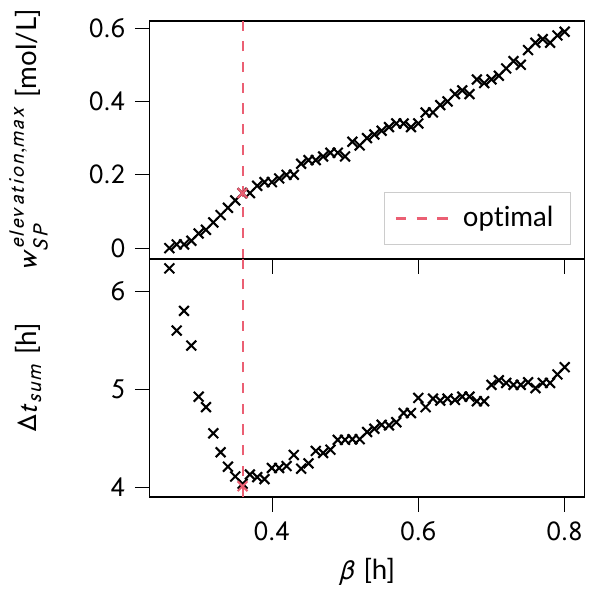}
    \caption{Result of the parameter tuning: maximum allowable set-point elevation $w_{SP}^{elevation,max}$ (top) and sum of all 6 transition times $\Delta t_{sum}$ resulting with $w_{SP}^{elevation,max}$ (bottom) for different values of the time constant $\beta$. We choose the optimal, i.e., smallest $\Delta t_{sum}$ by setting $\beta=0.36~h$ and $w_{SP}^{elevation}=0.15~\frac{mol}{L}$.}
   \label{fig:fit_beta}
\end{figure}

As Figure \ref{fig:fit_beta} shows, exploring the trade-off between set-point elevation and time constants improves the scale-bridging model performance significantly. The smallest, i.e., fastest, possible time constant $\beta^{min} = 0.26~h$, which  does not allow any set-point elevation, leads to a combined transition time of $\Delta t_{sum} = 6.23~h$.
The slightly higher time constant $\beta = 0.36~h$ in combination with a set-point elevation of $w_{SP}^{elevation} = 0.15~\frac{mol}{L}$ allows to reduce the transition time by 35\% to the optimum of $\Delta t_{sum} = 4.04~h$. 
The resulting transitions with the chosen optimal values are shown in Figure \ref{fig:transitions}. 
Note that the set-point elevation is not strictly increasing with $\beta$ due to the discretization. For example, with $\beta=0.49~h$, a set-point elevation of $w_{SP}^{elevation} = 0.26~\frac{mol}{L}$ is feasible, whereas for $\beta=0.50~h$, $w_{SP}^{elevation} = 0.26~\frac{mol}{L}$ is not feasible. 
The reason is that in one transition the set-point is at $w_{SP}^{max} = 0.76~\frac{mol}{L}$ for one discretization time step longer with $\beta=0.50~h$ compared to $\beta=0.49~h$ leading to a slight overshoot of the concentration out of the product band.

\subsubsection*{Energy demand model }
As discussed in the Method Section, 
model \circled{2} is needed in scheduling optimization to predict the process energy demand.
In this case study, model \circled{2} needs to predict the cooling power $Q_{cool}$ as a function of the controlled variable $y_{cv}$, which is  the concentration $C_A$, and its time derivatives. 
Note that, in principle, the cooling power $Q_{cool}$ is a process degree of freedom; however, the cooling power is set by the underlying PID-controller (\cref{eq:PID}).
In scheduling optimization, the degree of freedom  is the piece-wise constant set-point $w_{SP}$ for the concentration $C_A$.
This piece-wise constant set-point $w_{SP}$ is filtered, resulting in a reference trajectory for both $C_A$ and its time derivatives (compare to Figure~\ref{fig:method}).
Thus, the data-driven energy demand model must approximate the response of the closed-loop system as a function of $C_A$ and its time derivatives.

As the operation of the multi-product reactor is divided in production and transition periods and we need to model the cooling power accurately in particular during the long production periods, we split $Q_{cool}$ into a steady-state and a dynamic part:
\begin{align}
    \label{eq:q_cool_split}
    Q_{cool} = Q_{cool}^{steady} + Q_{cool}^{dynamic}
\end{align}
To approximate the first contribution $Q_{cool}^{steady}$, we assume that steady-state cooling powers are known for all three products and interpolate $Q_{cool}^{steady}$ as a piece-wise affine function of $C_A$. The three steady-state operating points $C_A = \{0.1~\frac{mol}{L}, 0.3~\frac{mol}{L}, 0.5~\frac{mol}{L}\}$ with corresponding cooling powers $Q_{cool}^{steady}$ lead to two piece-wise affine segments: The first affine segment approximates $Q_{cool}^{steady}$ for $C_A\leq 0.3~\frac{mol}{L}$ and the second affine segment approximates $Q_{cool}^{steady}$ for $C_A\geq 0.3~\frac{mol}{L}$. With the binary variable $z_{cool}^{steady}$, we can express $Q_{cool}^{steady}$ as
\begin{align}
    \label{eq:q_cool_stead}
    Q_{cool}^{steady} =& Q_{cool,0.3 \frac{mol}{L} }^{steady} + m_1^{steady}(1-z_{cool}^{steady})(C_A-0.3~\frac{mol}{L})   + m_2^{steady}z_{cool}^{steady}(C_A-0.3~\frac{mol}{L}),
\end{align}
where $Q_{cool,0.3\frac{mol}{L}}^{steady}$ is the steady-state cooling power at $C_A=0.3~\frac{mol}{L}$, $m_1^{steady}$ is the slope for $C_A\leq 0.3~\frac{mol}{L}$, and $m_2^{steady}$ is the slope for $C_A\geq 0.3~\frac{mol}{L}$. The slopes $m_1^{steady}$ and $m_2^{steady}$ are calculated from the cooling power at steady-state operating points (Table \ref{tab:products}). The bilinear terms $z_{cool}^{steady}C_A$ are reformulated using the Glover reformulation \citep{Glover.1975}.
                
\begin{figure}[h]
\centering
    \includegraphics{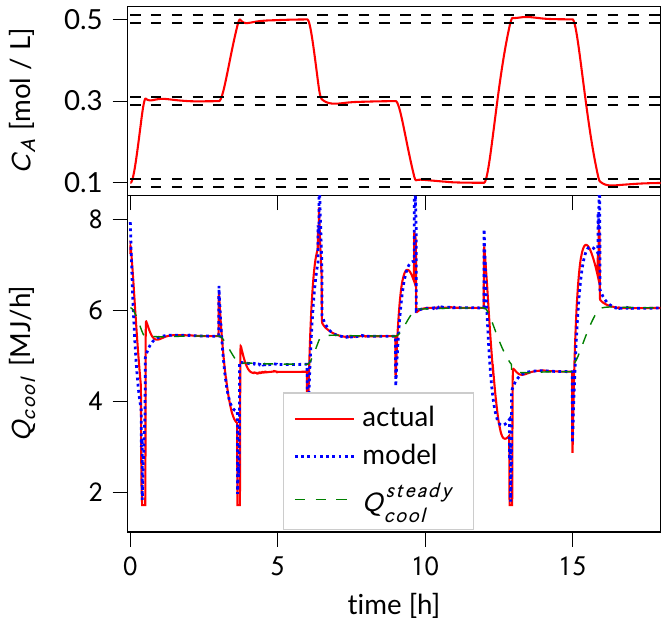}
    \caption{Fitting results for steady-state and total cooling power $Q_{cool}^{steady}$ and $Q_{cool}$, respectively, of the process energy demand model (compare to equations (27) and (28))}
   \label{fig:modelfit}
\end{figure}

The approximation of $Q_{cool}^{dynamic}$ is fitted to simulation data. Again, we simulate all six possible transitions using the nonlinear reactor model and the underlying control. The resulting cooling power is the red curve in Figure \ref{fig:modelfit}. 
The total cooling power deviates from $Q_{cool}^{steady}$ (dashed green curve in Figure \ref{fig:modelfit}) during transitions.
We model the dynamic part of the cooling power $Q_{cool}^{dynamic}$ as a linear function of the derivatives of the concentration, i.e.,
\begin{align}
    \label{eq:q_cool_dynamic}
    Q_{cool}^{dynamic} = c_1 \frac{dC_A}{dt} + c_2 \frac{d^2C_A}{dt^2},
\end{align}
with the two fitting parameters $c_1$, $c_2$, whose values are determined using the normal equation method \citep{Lewis.2006}.
The values are: $c_1 = -3.10~\frac{MJL}{mol}$, $c_2 = 0.444~\frac{MJLh}{mol}$. The resulting approximation of $Q_{cool}$ is shown in blue in Fig \ref{fig:modelfit}.
Note that in Figure \ref{fig:modelfit}, the concentration does not reach steady state after entering the product bands.
Still, the fitted dynamic cooling power $Q_{cool}^{dynamic}$ is negligible for 5 of 6 production phases and only in the second production phase a small offset between model and actual cooling power occurs. 
The purely linear model in \cref{eq:q_cool_dynamic} is the simplest possible choice to model the dynamic part and already leads to satisfying results.
In case study 3, a more complex model accounting for internal process states is demonstrated.

\subsubsection*{Energy system model}
For the energy system model \circled{3}, we have to calculate the electric input power $P_{CC,i}$ needed for the compression chillers as a function of the required cooling power $Q_{CC,i}$.
As the COP of the compression chillers depends on the part-load fraction (\cref{eq:case_study_3}), $P_{CC,i}$ is a nonlinear function of $Q_{CC,i}$, which we approximate as a piece-wise affine function.
We use two piece-wise affine segments per chiller with the breakpoint at 70~\% part load; two segments provide a good approximation \citep{Voll.2013}. 
The piece-wise affine curves can be modeled without introduction of additional binary variables, as the electric input power $P_{CC,i}$ is a convex function of the cooling power $Q_{CC,i}$.
Using equations from \cite{VerenaNeisen.2018}, we introduce two continuous variables $y_{i,1}, y_{i,2}$ that cover the two affine segments:
\begin{align}
    Q_{CC,i} &= y_{i,1} + y_{i,2}  &\forall i =1,2,3\\
    \label{eq:powercurve}
    P_{CC,i} &= P_{CC,i}^{min}z_{on,CC,i} + y_{i,1} m_{i,1} +  y_{i,2} m_{i,2} &\forall i =1,2,3
\end{align}
In \cref{eq:powercurve}, $P_{CC,i}^{min}$ is the electric input power at minimum part-load of chiller $i$, $z_{on,CC,i}$ is a binary indicating whether chiller $i$ is active, and $m_{i,1}, m_{i,2}$ are the slopes within the two piece-wise affine segments. 
Finally, we include the energy balance stating that the cooling demand of the reactor must be matched by the compression chillers:
\begin{align}
    \sum_{i=1}^{3} Q_{CC,i} = Q_{cool}
\end{align}

Further details on the scheduling optimization problem such as discretization and problem-specific constraints are given in the Supplementary Information (SI).

\subsection*{Sequential steady-state scheduling benchmark}
This benchmark represents a typical sequential scheduling without DR, referred to as SEQ in the following. 
First, the process schedule is optimized with only the product revenue $\Phi_{Product}$ in the objective function.
Second, the energy costs are minimized for fixed production decisions.
Detailed information are given in the SI.

\subsection*{Scheduling with full nonlinear model}
To estimate an upper bound on the economic performance of simultaneous scheduling, we perform an optimization with the nonlinear full-order system model.
To this end, we replace the models \circled{1}, \circled{2}, \circled{3} in the optimization problem with the nonlinear reactor model (\cref{eq:case_study_1,eq:case_study_2}) and the nonlinear compression chiller efficiency (\cref{eq:case_study_3}). Again, time is discretized using collocation and we receive a MINLP. 
We solve the MINLP optimization problem using BARON version 21.1.13 \citep{Khajavirad.2018} in heuristic mode, i.e., the resulting solution is no rigorous bound.
We refer to this benchmark as MINLP.
To obtain a feasible initial point, we fix the binary variables to the values resulting from our simultaneous dynamic scheduling and solve the resulting NLP.

\subsection*{Results}
\label{sec:results}
In this Section, we compare the economic profit obtained with our \emph{simultaneous dynamic scheduling} (SDS) to the sequential scheduling (SEQ) and the full-order nonlinear scheduling (MINLP).
While in case of the MINLP the profit is the objective value in the optimization, for the sequential scheduling and the \emph{simultaneous dynamic scheduling}, the profit is derived from a simulation of the original nonlinear process model. Accordingly, the optimized set-point sequence is used as input to a simulation of the underlying controller and the nonlinear process model.

\begin{figure}[h]
\centering
    \includegraphics{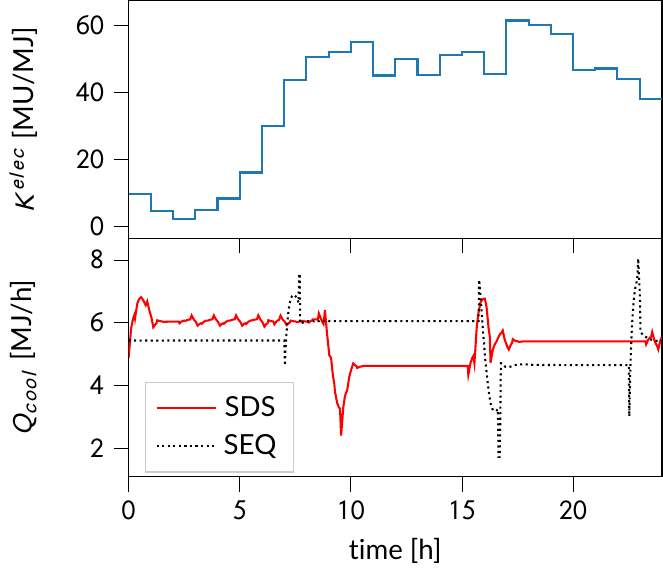}
    \caption{Electricity price $K^{elec}$ and simulated cooling power $Q_{cool}$ for \emph{simultaneous dynamic scheduling} (SDS) and sequential steady-state scheduling (SEQ). SDS performs  demand response and shifts cooling power to times of favorable prices.  }
   \label{fig:DR}
\end{figure}

The MINLP solution improves the SEQ solution by 5.5~\%. 
Our \emph{simultaneous dynamic scheduling} gains 5.2~\% compared to SEQ and thus captures 95~\% of the MINLP improvement. The improved economics mainly stem from demand response, i.e., products with higher cooling demands are produced at times of lower electricity prices (Figure \ref{fig:DR}). 
Additionally, we notice a higher energy efficiency during transition phases such that our 
\emph{simultaneous dynamic scheduling} reduces the amount of electricity consumed by 1.2~\% compared to SEQ.

\begin{figure*}[ht] 
\centering
    \includegraphics{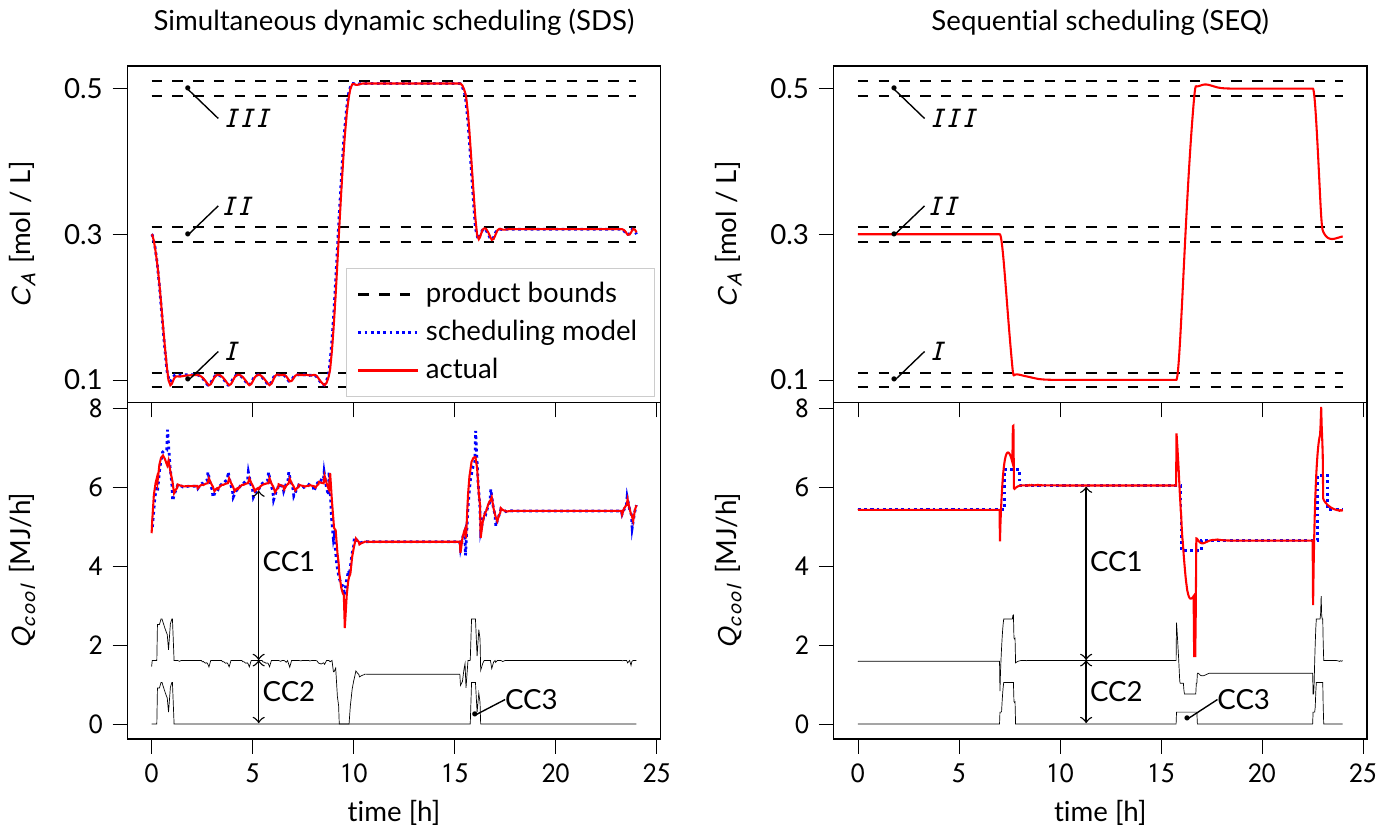}
    \caption{Comparison of concentration $C_A$ and cooling power $Q_{cool}$ between simultaneous dynamic scheduling (SDS, left) and sequential scheduling (SEQ, right). We indicate the three product bands ($\RNum{1}$, $\RNum{2}$, $\RNum{3}$) and the fraction of the cooling power $Q_{cool}$ that is supplied by the three compression chillers (CC1, CC2, CC3).}
   \label{fig:compCAQcool}
\end{figure*}

Figure \ref{fig:compCAQcool} shows concentration $C_A$ and cooling power $Q_{cool}$ for both SDS and SEQ.
In the case of SDS, the difference between the optimization model and the actual cooling power during the transitions is smaller than in the case of SEQ because the dynamics of the cooling power are modeled in SDS while SEQ only considers the average cooling power during a transition.
Modeling the dynamics of the cooling power within a transition leads to better scheduling decisions regarding the  on/off status of the three compression chillers.
The most distinct difference occurs in the transition from product $\RNum{1}$ to product $\RNum{3}$.
In the sequential scheduling, this transition features a high cooling power peak, which requires to turn on chiller 3 with the worst COP. 
In the case of SDS, the same transition is shaped such that it is not necessary to turn on chiller 3. Moreover, SDS anticipates that 
chiller 2 with the medium COP can be turned off during the second half of the transition. 
We expect that energy efficiency improvements are even higher in cases with longer transitions.
Naturally, the potential for energy efficiency improvements depends on the accuracy of the data-driven energy demand model (\circled{2}). 

Note that in our illustrative example, we assume that once the energy system components are active they can react instantaneously. Moreover, we assume that the frequency of on/off-switches resulting from the scheduling optimization with 15 minute resolution for discrete variables is acceptable. In practice, it might be necessary to consider ramp limits \citep{sass2019optimal}, or minimum up and down times \citep{Carrion.2006}. Such constraints can be added in a straight-forward manner to the formulation if needed.

\begin{figure}[h] 
\centering
    \includegraphics{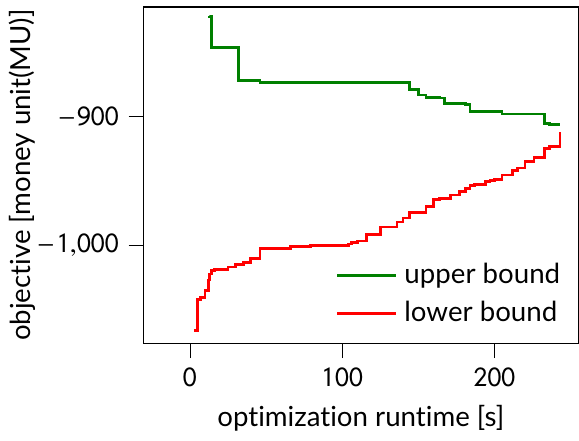}
    \caption{Convergence plot of \emph{simultaneous dynamic scheduling}.}
    \label{fig:convergence}
\end{figure}

A solution with a 1.0~\% optimality gap is found and proven in 244~s. 
Such a solution time is applicable for both offline day-ahead scheduling and online re-scheduling during the day, e.g., with a sampling time of one hour.
Note that in re-scheduling the solution from the last scheduling-iteration can be used for initialization to further speed up the optimization.
We also observe that SDS finds good feasible solutions quickly as shown in the convergence plot (Figure \ref{fig:convergence}).
After 32~s, a solution is found that has only 4.4~\% gap to the final lower bound and already outperforms the sequential scheduling. 
 
\subsubsection*{Influence of time constant $\beta$}

\begin{figure}[h] 
\centering
    \includegraphics{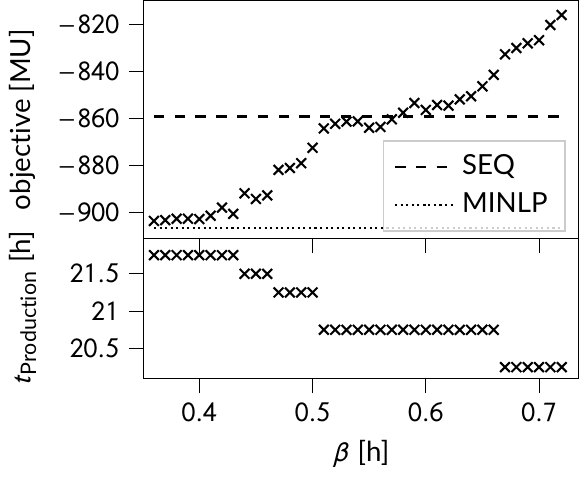}
    \caption{Objective (in money unit (MU)) and total production time $t_{\text{Production}}$ for different time constants $\beta$ starting from nominal $\beta = 0.36$~h up to $\beta = 0.72$~h. The objective values resulting from the sequential approach (SEQ) and MINLP optimization with both 21.75 hours of production are shown for comparison (dashed and dotted lines, respectively).}
    \label{fig:sens_beta}
\end{figure}

If a nonlinear process model is not available, the time constant $\beta$ cannot be chosen as in this case study 
but needs to be chosen based on intuition or recorded product transitions.
In result, the time constant might be suboptimal.
To study the influence of the time constant choice on the profit, the time constant in our case study is increased from the optimal value $\beta = 0.36$~h by up to 100~\% (Figure~\ref{fig:sens_beta}).
Note that the profit is calculated in a simulation using the original process model which leads to small differences between scheduling optimization and process simulation and therefore the simulated objective shown in Figure~\ref{fig:sens_beta} does not strictly increase with $\beta$. 

We find that as long as $\beta$ is increased by 20~\% or less the objective does not worsen more than 0.5~\%.
This result can be explained by the total production time which is 21.75~h for $\beta = 0.36$~h.
These 21.75 hours of production are still reached for $\beta = 0.43$~h and the loss in profit is small.
Generally, production time changes in 0.25~h steps corresponding to the time discretization of the binary variables (compare  to the SI).
When $\beta$ is increased by more than 20~\% above the optimal value, production is lost and the objective substantially worsens.
But, even for a 50~\% increase, the objective is still better than that of the sequential solution.

\section*{Case Study 2: Reactor with variable concentration}
While multi-product processes are one example for scheduling-relevant dynamics, single-product processes can also introduce dynamics if they can vary their controlled variable around a nominal value as long as the nominal value is reached on average over the considered time horizon.
To demonstrate that our SDS also works for a single-product case, we present a second case study and again study the influence of the time constant. 
The second case study is constructed by modifying the first one.

\subsection*{Setup}
A similar setup is used as in the previous case study with a CSTR and three compression chillers (Figure \ref{fig:case_study}).
Instead of a multi-product CSTR, we assume a single-product CSTR with a nominal concentration $C_A^{nom} = 0.3~\frac{mol}{L}$.
The CSTR has flexibility because we assume that the concentration is allowed to vary between $C_A^{min} = 0.09~\frac{mol}{L}$ and $C_A^{max} = 0.51~\frac{mol}{L}$ as long as the nominal concentration is reached on average over the time horizon $t_f - t_0$.
Accordingly, the condition 
    \begin{align}
        \label{eq:condition_integral}
        \int_{t_0}^{t_f}C_A  dt = C_A^{nom}(t_f-t_0)
    \end{align}
    has to hold.
Such setups occur when the product can be stored and is well-mixed in the storage tank.  
Note that an equation of this type would also occur for processes having a variable production rate as controlled variable $y_{cv}$. Only the concentration $C_A$ in \cref{eq:condition_integral} would have to be replaced by the production rate.

Obviously, it is favorable to operate at concentrations with a high cooling demand at times of low electricity prices and at concentrations with low cooling demand at times of high prices.
The challenge for scheduling optimization thus is to find a trajectory for the concentration that (i) can be realized by the process and (ii) reaches the nominal concentration on average. At the same time, the on/off status of the chillers has to be determined.
Note that the scheduling problem has three differences compared to the previous case:
    \begin{enumerate}
        \item Only the cumulative energy costs at final time $\Phi_{Energy}(t_{f})$ are considered in the objective (\cref{eq:obj}) since the production volume is fixed.
        \item All constraints associated with the different products and the production bands are removed (equations (S1), (S2), (S6)-(S10) in the SI)
        \item \Cref{eq:condition_integral} is included such that the nominal concentration is reached on average.
    \end{enumerate}
As energy costs are the only objective function in the second case study, a sequential scheduling is not applicable because there is no objective for the process optimization.
Thus, we benchmark our SDS against a steady-state operation of the CSTR at the nominal concentration. Again, as a second benchmark, a MINLP optimization is performed using BARON in heuristic mode.
A feasible initial point is found by fixing the binary variables to the values from our \emph{simultaneous dynamic scheduling} and solving the resulting NLP.

\subsection*{Results}

The MINLP solution improves the steady-state solution by 6.7~\%. Our \emph{simultaneous dynamic scheduling} (SDS) reduces costs by 5.5~\% compared to the steady-state benchmark and thus captures 82~\% of the MINLP improvement. 
The optimization runtime of our SDS approach is only 55~s.
Note that again the MINLP optimization with BARON does not provide a feasible point without initialization from the SDS solution.

\begin{figure*}[h] 
\centering
    \includegraphics{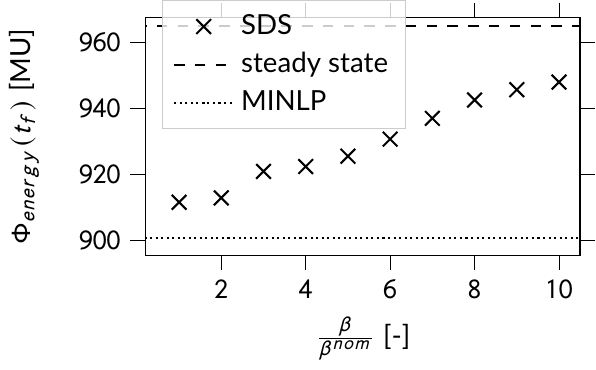}
    \caption{Energy costs $ \Phi_{energy}(t_f)$ (in money unit (MU)) in the second case study achieved with simultaneous dynamic scheduling (SDS) for different time constants $\beta$ normalized to nominal value $\beta^{nom}=0.36$~h. The energy costs resulting from steady-state operation and the MINLP benchmark are shown for comparison (dashed and dotted lines, respectively).}
    \label{fig:sens_beta_no_product}
\end{figure*}

\begin{figure*}[h] 
\centering
    \includegraphics{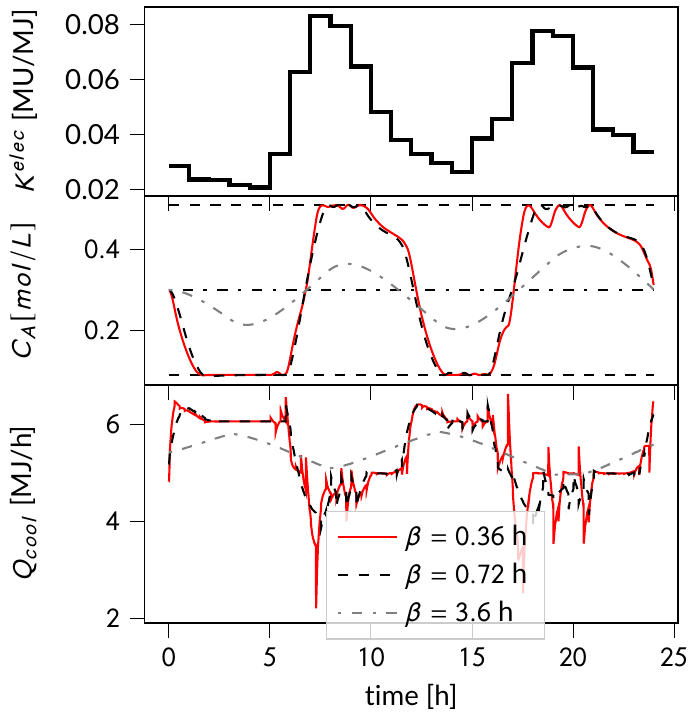}
    \caption{Electricity price $K^{elec}$, concentration $C_A$, and cooling power $Q_{cool}$ in the second case study for three values of the time constant $\beta$.}
   \label{fig:operation_no_product}
\end{figure*}

Compared to case study 1, the choice of the time constant $\beta$ has a much lower impact on the economic result (Figure \ref{fig:sens_beta_no_product}).
Note that we use the same time constant $\beta$ as in case study one, because the transitions studied during tuning cover the complete range of allowed concentrations.  
If $\beta$ is doubled from 0.36~h to 0.72~h, the cost reduction still amounts to 5.4~\% compared to steady-state operation (Figure \ref{fig:sens_beta_no_product}).
The operation is very similar for both time constants and the cooling power only deviates significantly in hours 1, 6-7, 11, and 16-17 (Figure \ref{fig:operation_no_product}).
In hours 16 and 17, the schedule with $\beta = 0.72$~h drives the reactor from minimum concentration to maximum concentration. 
The higher flexibility of the low time constant $\beta = 0.36$~h allows to consume more cooling in hour 16 and less in hour 17 compared to the case where $\beta = 0.72$~h.
Still, the larger time constant can capture the main trend of the electricity price profile, which has a peak in the morning and another one in the afternoon.
Such a price profile is typical for the German market where the main price periodicities are 24 and 12 hours \citep{Schafer.2020}.
Even if $\beta$ is increased by a factor of 10, the scheduling can still capture the main trend of the electricity price profile (Figure \ref{fig:operation_no_product}).
Therefore, significant cost reductions can be reached even if the chosen time constants are far above the optimal value.

\section*{Case Study 3: Distillation column with variable purity}
In this Section, we study a case similar to the second case  but, instead of a SISO CSTR, we consider a distillation column as a MIMO process. 
The purity of both top and bottom product can be varied throughout the day as long as the desired purity is reached on average. 

\subsection*{Setup}

\begin{figure}[h]
\centering
    \includegraphics{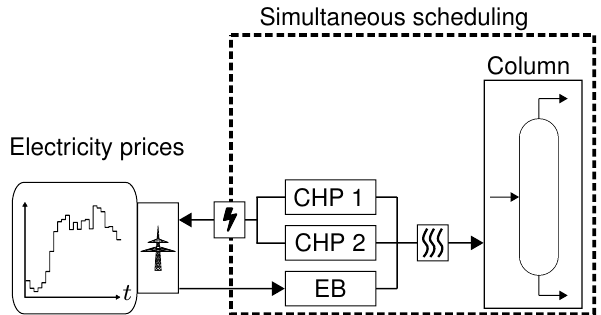}
    \caption{Case study 3: Simultaneous scheduling of a distillation column heated by two combined heat and power plants (CHP1, CHP2) and a electricity-driven boiler (EB). Time-varying electricity prices provide an economic incentive for DR. }
    \label{fig:case_study_col}
\end{figure}

The heat demand of a distillation column is satisfied by two combined heat and power plants (CHPs) and an electricity-driven boiler (EB) (Figure \ref{fig:case_study_col}).
Electricity produced by the CHPs is sold to the electricity grid.
For the distillation column, we use a generic benchmark model of a binary distillation proposed by  \cite{Skogestad.1988} together with a liquid flow model from \cite{Skogestad.1990}. 
The column  model consists of $N+1$ mass balances of the light boiling component and $N+1$ mass balances of the liquid hold-ups where $N=40$ is the number of theoretical trays.
In total, the column model has 82 differential states.
It is assumed that the heat demand of the column is proportional to the boilup flow rate $V$ and the heat demand is scaled such that the column requires 1 MW of heating in nominal operation.
The two CHPs have 800~kW and 500~kW nominal thermal power and are subject to 50~\% minimum part-load constraints \citep{Voll.2013}.
The electricity-driven boiler has 800 kW nominal power and a 20~\% minimum part-load constraint \citep{Baumgartner.2020}.
Further details on the distillation model and the CHP and EB models are given in the Supplementary Information (SI).

For the column, the feed flow $F$ is fixed by an upstream process.
Four flows can be manipulated: the reflux flow rate $L$, the boilup flow rate $V$, the distillate flow rate $D$, and the bottom flow rate $B$. 
Accordingly, the vector of manipulated variables is 
\begin{align}
    \vect{u} = \left( \begin{array}{c} L \\ V \\ D \\ B \end{array} \right).
\end{align}
The variables to be controlled are the vapor mole fraction of the light component entering the condenser $y_D$, its liquid mole fraction in the bottom flow $x_B$, the condenser hold-up $M_D$, and the condenser hold-up $M_B$:
\begin{align}
    \vect{y}_{cv} = \left( \begin{array}{c} y_D \\ x_B \\ M_D \\ M_B \end{array} \right)
\end{align}
We couple the mole fractions $y_D$ and $x_B$ with each other by defining the purity $\rho$ which is equal to $y_D$  and $1-x_B$.
The coupling $y_D = 1-x_B$ can be applied in this case study, as the feed purity $z_F$ is 50~\% and in steady-state both bottom flow rate $B$ and the distillate flow rate $D$ are equal to 50~\% of the feed flow rate $F$. 
A more general case, is discussed in the SI.

The hold-ups $M_D$ and $M_B$ shall be maintained constant at their nominal values $M_D^{nom}$ and $M_B^{nom}$ irrespective of the current purity $\rho$.
Accordingly, the vector of filtered set-points $\vect{w}_{SP,fil}$ can be given as a function of the filtered purity set-point $\rho_{SP,fil}$:
\begin{align}
    \vect{w}_{SP,fil}(\rho_{SP,fil}) = \left( \begin{array}{c} \rho_{SP,fil}\\ 1-\rho_{SP,fil} \\ M_D^{nom} \\ M_B^{nom} \end{array} \right)
\end{align}
The purity $\rho$ can be varied between $\rho^{min}=0.85$ and $\rho^{max}=0.95$ as long as the nominal value $\rho^{nom}=0.9$ is reached on average.

To this end, we use two PI controllers \citep{Corriou.2018} to control the mole fractions $y_D$ and $x_B$ by manipulating the flows rates $L$ and $V$, respectively.
Further details are given in the SI.
For the hold-ups $M_D^{nom}$ and $M_B^{nom}$, we do not model the controllers explicitly but follow the common assumption that an underlying control sets the flows $D$ and $B$ such that the hold-ups are controlled perfectly \citep{Skogestad.1988}.

With this case study, we demonstrate that a MIMO process does not necessarily lead to a MIMO scale-bridging model (cf., discussion in the Method Section). 
Instead, controlled variables such as the mole fractions $y_D$ and $x_B$ can be given as functions of a single scheduling-relevant variable such as the purity $\rho$ here.
At the same time, other controlled variables, such as the hold-ups $M_D$, $M_B$, need to be maintained constant irrespective of the scheduling-relevant variable.
As a consequence, instead of four SBMs only one SBM describing the dynamics of the purity $\rho$ is needed.

\subsection*{Simultaneous dynamic scheduling}
In this Section, we develop the three parts of our model.

For the scale-bridging production process model \circled{1}, a first-order model is chosen because all 4 controlled variables in $\vect{y}_{cv}$ can be controlled with a relative degree $r=1$ (cf. differential model equations in the SI).
Thus, model \circled{1} for the filtered  purity set-point $\rho_{SP,fil}$ is
\begin{align}
   \rho_{SP,fil} + \tau \frac{d \rho_{SP,fil}}{dt} = \rho_{SP},
\end{align}
with one time constant $\tau$ that must be tuned. 
For the tuning, a similar procedure as in the previous case studies is applied. 
Again, we study 6 representative transitions.
The aim is to find a time constant $\tau$ and set-point bounds $\rho_{SP}^{min}$, $\rho_{SP}^{max}$ that give fast transitions.
At the same time, deviations between $y_D$ and $\rho_{SP,fil}$ and between $x_B$ and $1-\rho_{SP,fil}$ need to be within a certain tolerance.
Details on the tuning are given in the SI. 
The resulting values are $\tau= 10$~min, $\rho_{SP}^{min} = 0.8485$, $\rho_{SP}^{max} = 0.9515$.

\begin{figure}
    \centering
    \includegraphics{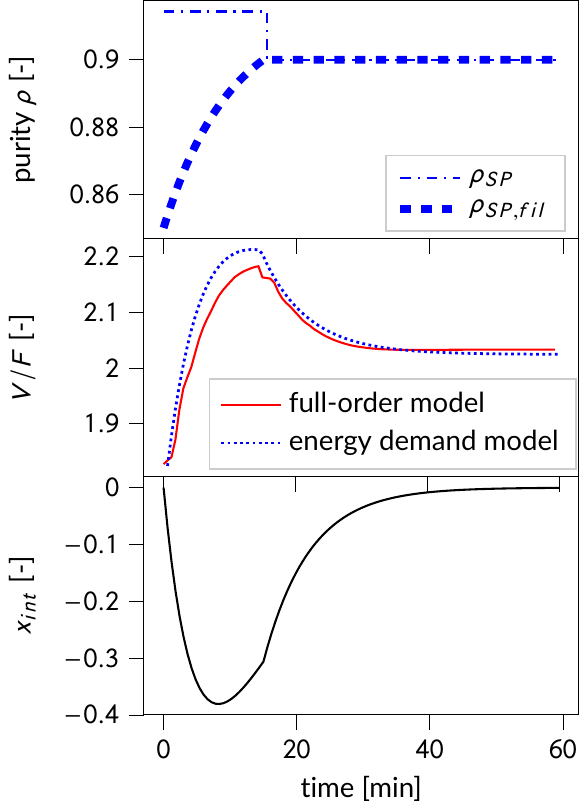}
    \caption{Exemplary transition showing purity set-point $\rho_{SP}$ and filtered set-point $\rho_{SP,fil}$ (top), reboiler flow $V$ normalized to the feed flow $F$ calculated from full-order model and energy demand model (\cref{eq:V_split,eq:V_dynamic}) (middle), and the internal state $x_{int}$ of the energy demand the model calculated based on \cref{eq:dx_int_dt} (bottom).}
    \label{fig:fit_x_int}
\end{figure}

For the energy demand model \circled{2}, we split the vapor flow $V$ into a steady-state and a dynamic part (compare to \cref{eq:q_cool_split}), i.e.,
\begin{align}
    \label{eq:V_split}
    V = V^{steady} + V^{dynamic},
\end{align}
where $V^{steady}$ is modeled as a piece-wise affine function with two segments of the purity $\rho$ because the vapor flow $V$ in steady-state is a nonlinear function of the purity.
The dynamic part $V^{dynamic}$ is approximated by a linear model with the  derivative of the scale-bridging variable, $\frac{d\rho_{SP,fil}}{dt}$, as input.
However, in contrast to \cref{eq:q_cool_dynamic}, the model for $V^{dynamic}$ features an internal state $x_{int}$ and is given by
\begin{align}
    \label{eq:dx_int_dt}
    &\frac{dx_{int}}{dt} = a x_{int} + b\frac{d\rho_{SP,fil}}{dt},\\
     \label{eq:V_dynamic}
    &V^{dynamic} = cx_{int} + d\frac{d\rho_{SP,fil}}{dt},
\end{align}
with fitting coefficients $a, b, c, d$ that are fitted to the simulated transitions that result from the tuning of model \circled{1} discussed above.
The rational for including an internal state $x_{int}$ is shown in Figure~\ref{fig:fit_x_int} where one transition is visualized:
After 15 minutes, the filtered purity set point $\rho_{SP,fil}$ has reached the steady state value, i.e., $\rho_{SP}  = \rho_{SP,fil}$ and $\frac{d\rho_{SP,fil}}{dt} = 0$.
However, the vapor flow $V$ is not in steady state because the uncontrolled states inside the column are not in steady state and it takes roughly another 15 minutes until steady state is reached. 
The fitting constant $a = -0.14~\frac{1}{min}$  corresponds to a time constant of 7.1 minutes and represents internal dynamics of the column. 
We also studied energy demand models with more than one internal state but did not find significant improvements compared to the one with a single internal state (\cref{eq:dx_int_dt}).

With this third case study, we demonstrate the model-order reduction capabilities of scale-bridging models and data-driven models. 
The 82 differential states of the full-order column model are reduced to just 2 states: the filtered purity set-point $p_{SP,fil}$ and the internal state of the energy demand model $x_{int}$.

In the energy system model \circled{3}, the part-load behavior of the CHPs is modeled with one piece-wise affine segment leading to a reasonable discretization \citep{Voll.2013}.
For the electricity-driven boiler, a constant efficiency is assumed \citep{NilsBaumgartner.2019}.

We use a quarter-hourly electricity price profile because $\Delta t_{elec}= 15~min$ being similar to the time constant of the SBM ($\tau = 10~min$) presumably causes more pronounced dynamic operation than in the previous case studies. 
The price profile is taken from the German intra-day market and occurred on the 13th January, 13th, 2021 \citep{energyCharts}. 
A discretization with $\Delta t_{cont}= 15~min$ and $N_{cp} = 3$ collocation points is used.
As the 15-minute price profile causes excessive on-off switching of the energy system components in preliminary investigations, we require a minimum up-time of one hour and a minimum down-time of one hour for all three energy system components.
Such constraints are often used in energy system optimization as excessive on-off switches often have negative effects on component life-time \citep{Carrion.2006}.
The scheduling optimization problem is solved with gurobi 9.1.2 \citep{gurobi} on the same machine as before.
Again we use a 1~\% optimality gap.

\subsection*{Benchmarks}
We consider two benchmarks, a steady-state operation at nominal purity and a sensitivity analysis where we halve the nominal time constant $\tau = 10~min$.
The latter leads to operations with time constants that render some transitions infeasible but allows us to estimate the performance lost by the scale-bridging approach in comparison to an optimization with the original full-order nonlinear process model that, unlike the SBM, would not be limited by the slowest transition.
We perform this sensitivity analysis with smaller time constants since an optimization with the full-order process model in this case study would be too computationally challenging due to the large scale of the MINLP.

\subsection*{Results}
Our simultaneous dynamic scheduling (SDS) reduces costs by 4.3~\% compared to a steady-state operation whose results are shown in as Figure~S3 of the SI. 
The SDS optimization converges to the pre-defined optimality gap of 1~\% within 125 s. 
Like in the previous case studies, near-optimal feasible points are found even faster (cf. Figure~S4 in the SI).

\begin{figure}
\centering
    \includegraphics[width=15cm]{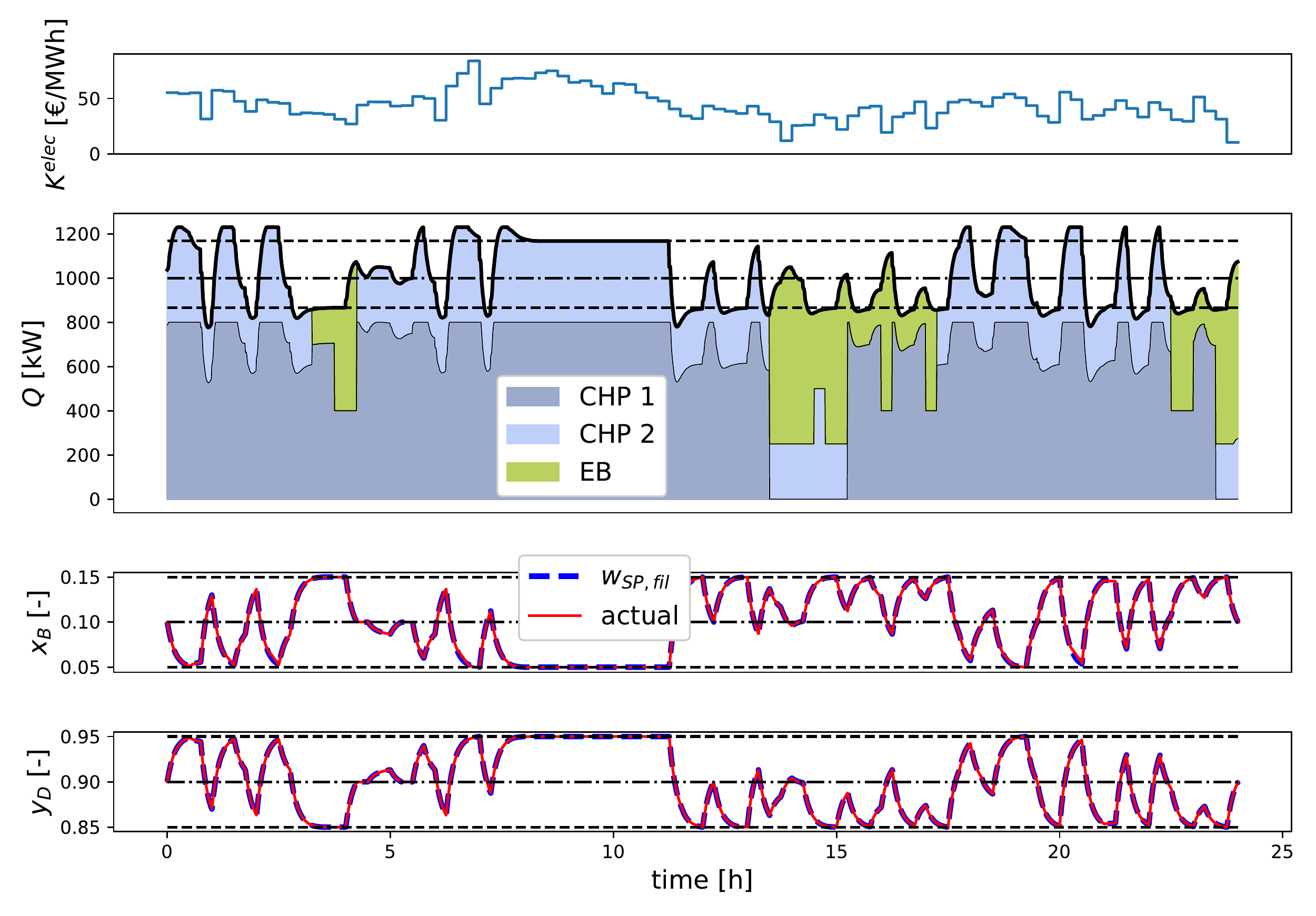}
    \caption{Resulting operation for the third case study. Top: Electricity price $K^{elec}$. Center: Heat demand $Q$ of the column (bold black line) together with nominal value (dashed dotted line) and minimum and maximum steady-state heat demands (dashed lines). The portions of the heat demand supplied by the two combined heat and power plants (CHP1 and CHP2) and the electrically-driven boiler (EB) are indicated with colors. Bottom: Actual values of bottom  composition $x_B$ and top composition $y_D$ together with their respective filtered set-points $w_{SP,fil}$.}
    \label{fig:case_study_3_operation}
\end{figure}

The resulting operation with the nominal time constant $\tau = 10~min$ is shown in Figure~\ref{fig:case_study_3_operation}.
Due to the 15-minute price profile,  the operation is more dynamic than in the previous case studies.
The column is only operated in steady state between $t = 3.5~h$ and $t = 4~h$ and between $t = 8.25~h$ and $t = 11.25~h$ .
Still, the PI-controllers track the filtered set-point accurately for both mole fractions $x_B$ and $y_D$ (Figure~\ref{fig:case_study_3_operation}).
The average mole fractions in the storage after the one day scheduling time horizons are:
\begin{align}
    &x_D^{average} = \frac{\int_{t = 0 h }^{t = 24 h} x_D D dt}{\int_{t= 0 h }^{t = 24 h} D dt} = 0.90040 > \rho^{nom}=0.9\\
    &x_B^{average} = \frac{\int_{t = 0 h }^{t = 24 h} x_B B dt}{\int_{t= 0 h }^{t = 24 h} B dt} = 0.09967 < 1-\rho^{nom}=0.1
\end{align} 
Thus,  both top and bottom product stream achieve the desired purity on average.

Figure~\ref{fig:case_study_3_operation} shows that the optimization seeks to operate the column at high purities and high heat demands at times of high prices (e.g., hour 20 - 20.5).
By supplying  the required heat from the CHPs, electric power can be sold to the market at a high price. 
However, also at times of low electricity prices, a comparatively high heat demand can be seen (e.g., hour 13.75 - 14.25 or 16 - 16.25).
Here, the electricity-driven boiler is operated close to its maximum load because  electricity is cheap. 
Only at times of medium prices, the column is operated at low purities resulting in a low heat demand (e.g. hour 3.25 - 4).
These interactions between the purity of the column and the on/off-status of energy system components demonstrate the advantages of a simultaneous scheduling of process and energy system.

\begin{figure}
\centering
    \includegraphics{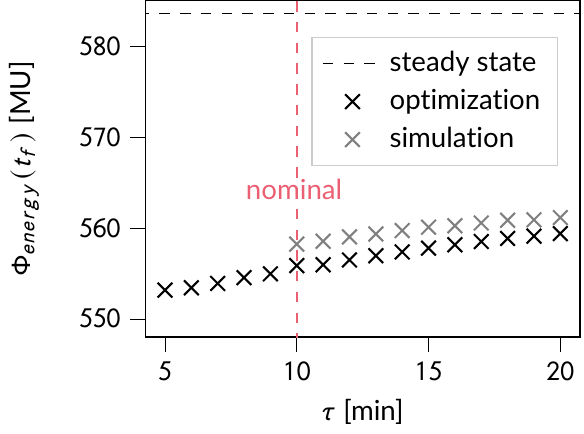}
    \caption{Energy costs $ \Phi_{energy}(t_f)$ (in money unit (MU)) in the third case study achieved with different time constants $\tau$. The energy costs resulting from steady-state operation are shown for comparison.}
    \label{fig:case_study_3_sensitivity}
\end{figure}

In Figure~\ref{fig:case_study_3_sensitivity}, the sensitivity of the energy costs with respect to the time constant $\tau$ is shown. 
For time constants above the nominal value, the optimized energy costs are compared to simulated energy costs.
The energy costs obtained in the simulation are always slightly higher than those predicted by the optimization due to inaccuracies of the data-driven energy demand model \circled{2}.
Comparing different time constants, we can conclude that the sensitivity towards the time constant is relatively weak. 
For instance, when the time constant is doubled from $\tau = 10~min$ to $\tau = 20~min$, the simulated energy costs only increase by 0.5~\%.

To estimate the performance lost due to the conservatism of the scale-bridging approach, we halve the nominal time constant $\tau = 10~min$ and study the resulting optimized energy costs.
Note that this yields a lower bound for the energy costs as some transitions are not feasible for the halved time constant $\tau = 5~min$.
However, this lower bound is only 0.5~\% smaller than the optimized energy costs with the nominal time constant (Figure~\ref{fig:case_study_3_sensitivity}).
Thus, the 4.3~\% cost reduction achieved by our SDS approach already realizes most of the overall demand response potential while being able to run within 5 minutes.
We therefore conclude that SDS offers a favorable compromise between solution quality and optimization runtime.

\section*{Conclusion and Discussion}
For power-intensive processes, volatile electricity prices provide an opportunity to increase profit via demand response (DR).
A particularly promising DR option is the simultaneous scheduling optimization of processes and their energy systems.
As such an optimization must consider scheduling-relevant process dynamics as well as on/off-decisions in the energy supply system, computationally challenging nonlinear mixed-integer dynamic optimization (MIDO) problems arise.
In this work, we present an efficient \emph{simultaneous dynamic scheduling} (SDS) approach that relies on a tailored scheduling model consisting of (i) a linear scale-bridging model for the closed-loop response of the process, (ii) a data-driven model for the process energy demand, and (iii) a mixed-integer linear programming (MILP) model of the energy system. 
Using a discrete time formulation and collocation, we receive an overall MILP formulation that can be optimized in practically relevant times.

First, we apply the method to a case study of a multi-product continuous stirred tank reactor (CSTR) cooled by three compression chillers.
Compared to a typical sequential scheduling, we find that the presented SDS
approach improves economic profit by 5.2~\%, just shy of the 5.5~\% found by nonlinear scheduling optimization using the original nonlinear process model. 
Second, we investigate a single-product reactor with a variable concentration. 
Here, SDS outperforms a steady-state operation by 5.5~\% while a nonlinear scheduling reaches 6.7~\%. 
Third, we investigate a distillation column heated by two combined heat and power plants and an electrically-driven boiler. 
The distillation column is a $4\times4$ multi-input multi-output (MIMO) process; however, one scale-bridging model is sufficient as we couple top and bottom purity and hold condenser and reboiler hold-up constant irrespective of the purity.
We thereby demonstrate the model-order reduction potential of our approach driving down the number of states from 82 to 2. Cost reduction of 4.3~\% compared to steady-state operation are achieved.

In all three case studies, the optimization runtime is sufficiently fast for online optimization.
As the proposed scheduling model always has the same basic structure, 
we expect the method to be real-time applicable in many cases.

A restriction of our method is that the scale-bridging approach imposes a single common linear closed-loop response in all operating regimes, which may cut off some of the process flexibility and thus DR potential.
For example, in our first case study, we must choose the time constants of the enforced linear closed-loop response such that all six product transitions are feasible.
Due to the nonlinear behavior of the CSTR, some of the transitions could in principle be performed faster, however, the critical transition, i.e., the slowest one, limits the time constants for the scale-bridging model. 
Moreover, it may in general be difficult to find the time constants that give the fastest possible linear closed-loop response. 
Finding the time constants using the heuristic used in our case studies is straightforward if the controlled process can be simulated and the relevant transitions can be studied in numerical experiments. 
If several controlled variables should be varied independently of each other, the combinatorial  complexity of the heuristic procedure would increase as more than one scale-bridging model would be needed.
Then, it might no longer  be possible  to simply sample the space of possible combinations but some kind of black-box optimization may be required.
For multi-product processes, which are inherently dynamic, the time constants can also be chosen based on recorded transitions. 
Our sensitivity study shows that, as long as transition times are only moderately larger than necessary, costs can still be reduced compared to a standard sequential scheduling. 
For processes that are currently operated in steady state without DR, no recorded transitions might be available. 
However, we demonstrate that, for such processes, time constant choice is less critical as even greatly suboptimal values may allow to follow slow trends in the electricity price profile.

Overall, our results demonstrate that the proposed method offers a favorable trade-off between accurate handling of dynamic flexibility and online applicable optimization run-times.
  
  \section*{Author contributions}
\textbf{Florian J. Baader}: Conceptualization, Methodology, Software, Investigation, Validation, Visualization, Writing - original draft. 
\textbf{André Bardow}: Funding acquisition, Conceptualization, Supervision, Writing – review \& editing. 
\textbf{Manuel Dahmen}: Conceptualization, Supervision, Writing – review \& editing.

  \section*{Declaration of Competing Interest}
    We have no conflict of interest.

    \section*{Acknowledgements}
    This work was supported by the Helmholtz Association under the Joint Initiative ``Energy System 2050--A Contribution of the Research Field Energy'' and under the Initiative ``Energy System Integration''.

\twocolumn
\section*{Nomenclature}
\label{sec:nomenclature}

\noindent\textbf{Abbreviations} \\
\noindent
\begin{tabularx}{\columnwidth}{lX}
    CC & compression chiller \\
    CHP & combined heat and power plant \\
    COP & coefficient of performance \\
    CSTR & continuous-stirred-tank-reactor \\
    DR & demand response \\
    EB & electrically-driven boiler\\
    KKT & Karush-Kuhn-Tucker \\
    MIDO & mixed-integer dynamic optimization \\
    MILP & mixed-integer linear programming \\
    MIMO & multiple-input multiple output \\
    MINLP & mixed-integer nonlinear programming \\
    MPC & model predictive control \\
    MU & money unit \\
    PID & proportional–integral–derivative \\
    PI & proportional–integral \\
    SBM & scale-bridging model \\
    SDS & simultaneous dynamic scheduling \\
    SEQ & sequential scheduling benchmark \\
    SI & Supplementary information \\
    SISO & single-input single output \\
    SO-MPC & scheduling-oriented model predictive control
\end{tabularx} \\

\noindent\textbf{Greek symbols} \\
\noindent
\begin{tabularx}{\columnwidth}{lX}
    $\beta$ &time constant of natural oscillation \\
    $\varepsilon$ & safety margin \\
    $\zeta$&damping coefficient \\
    $\varrho$ & density\\
    $\rho$ & scheduling-relevant variable \\
    $\tau$ & time constant \\
    $\overline{\tau}$ & scaled time \\
    $\Phi$ & objective 
\end{tabularx} \\

\noindent\textbf{Latin symbols} \\
\noindent
\begin{tabularx}{\columnwidth}{lX}
    $\vect{A}-\vect{H}$& matrices \\
    $a, b, c, d$ & fitting coefficients \\
    $B$ & bottom flow rate \\
    $C_A$ & concentration of component A \\
    $c_P$& heat capacity \\
    $D$ & distillate flow rate \\
    $E_A$ & activation energy \\
    $F$ & feed flow rate \\
    $f$, $g$& functions \\
    $f_e$ & finite element \\
    $\Delta H_r$ & enthalpy of reaction \\
    $K$&price \\
    $K_P$ & proportional controller constant \\
    $k$ & reaction constant \\
    $L$ & reflux flow rate \\
    $l$ & Lagrange polynomial \\
    $M_k$ & hold-up of tray $k$ \\
    \end{tabularx} \\
\begin{tabularx}{\columnwidth}{lX}
    $m$ & linear slope \\
    $N$ & number of trays \\

    $N_{cp}$ & order of collocation polynomial \\
    $n$ & natural number \\
    $P$ &power\\
    $Q$ & thermal power \\
    $q$ & flow rate \\
     $R$ &gas constant \\
    $r$ & order of differential equation\\
    $T$ & temperature \\
        $t$ & time \\
    $u$ & manipulated variable \\
    $V$ & volume (case study 1\&2) or boilup flow rate (case study 3) \\
    $w_{SP}$ & set-point\\

    $x$&differential state  \\
    $x_k$ & liquid mole fraction on tray k \\
    $y$ & continuous variable \\

    $y_k$ & vapor mole fraction on tray k \\
    $z$ & discrete variable
\end{tabularx} \\

\noindent\textbf{Sets} \\
\noindent
\begin{tabularx}{\columnwidth}{lX}
    $\mathbb{E}$ & end-energy forms \\
    $\mathbb{P}$ & products \\
    $\mathbb{S}$ &set of possible transitions \\
    $\mathbb{T}_{dis}$ & timepoints on discrete grid
\end{tabularx} \\

\noindent\textbf{Subscripts} \\
\noindent
\begin{tabularx}{\columnwidth}{lX}
    $0$ & initial\\
    $B$ & reboiler \\
    $c$ & component \\
    $cv$ & controlled variable\\
    $cont$ & continuous \\
    $cool$ & cooling \\
    $D$ & condenser \\
    $dis$ & discrete \\
    $e$ & end-energy form \\
    $ec$ &energy costs\\

    $ed$ & energy demand \\
    $elec$ &electricity \\
    $f$ & final \\
    $f_e$ & finite element \\
\end{tabularx} \\
\begin{tabularx}{\columnwidth}{lX}
    $fil$ & filtered\\
    $I$ & integral \\
    $in$ & input\\
    $int$ & internal \\
    $on$ & on-off status \\
    $out$ & output \\
\end{tabularx} \\
\noindent
\begin{tabularx}{\columnwidth}{lX}
    $p$ & product \\
     $s$ & transition  \\
     $SP$ & set point \\
   $sum$ & summed value
\end{tabularx} \\

\newpage
\noindent\textbf{Superscripts} \\
\noindent
\begin{tabularx}{\columnwidth}{lX}
    $end$ & final value \\
    $elec$ & electricity \\
    $max$ & maximum value\\
    $min$ & minimum value\\
    $nom$ & nominal \\
    $start$ & starting value \\
    $steady$ & steady-state value
\end{tabularx} \\
\onecolumn

\bibliographystyle{apalike}
\renewcommand{\refname}{Bibliography} 
\bibliography{literature.bib}

\end{document}


\thispagestyle{firststyle}

  \begin{center}
    \begin{large}
    {\fontsize{12}{14} \selectfont
      \textbf{\mytitle}}
    \end{large} \\
    \myauthor
  \end{center}

  \begin{footnotesize}
    \affil
  \end{footnotesize}

\section*{Case study 1: Multi-product reactor}

\subsection*{Details on the scheduling optimization problem}
\label{sec_opt_problem}
In this case study, we use the objective function as introduced in the Method Section of the main paper (equations (4) - (6)).
For time discretization, we have $\Delta t_{elec} = 1~h$ as we consider hourly changing electricity prices. 
For discrete variables, we choose $\Delta t_{dis} = 15$~min so the optimization can differentiate between transitions that can be performed within 30, 45, or 60 minutes. 
This resolution is reasonable as the transition times are between 27 and 54 minutes (Figure 6 of the main paper).
Note that though a finer discretization, e.g., 10 minutes, may lead to even higher profits, more binary variables will further increase the computational burden.
For the continuous variables, we find a discretization of $\Delta t_{cont} = 15$~min with $N_{cp}=3$ collocation points to give a sufficiently accurate time discretization.
In the main paper, we evaluate the optimization result in a simulation on the original nonlinear process model and thereby also verify the accuracy of our time discretization. 

\begin{table}[h]
    \centering
    \caption{CSTR model parameters adapted from \cite{Petersen.2017} (compare to equations (11) and (12) in the main paper)}
    \large
    \renewcommand\theadfont{\large}
    \label{tab:cstr_parameter}
    \begin{tabular}{lcc}
    \hline
    \thead{symbol} & \thead{value} & \thead{unit} \\
    \hline
    $q$ & $100$  & $\frac{m^3}{h}$\\
    $V$ &$100$  &$m^3$\\
    $C_{A,feed}$ &$1$  &$\frac{mol}{L}$\\
    $k$ & $7.2 \cdot 10^{10}$ & $\frac{1}{h}$\\
    $\frac{E_A}{R}$ & $6500$ & $K$\\
    $T_{feed}$ &$350$  &$K$\\
    $\frac{\Delta H_r}{R\varrho C_P}$ & $-209$  & $\frac{K m^3}{mol}$\\
    $\varrho$ &$1000$  &$\frac{kg}{m^3}$\\
    $C_P$ &$0.239$  &$\frac{J}{kgK}$\\
    \hline  
    \end{tabular}
\end{table}

Finally, we add problem-specific constraints. Two constraints per product determine if a product $p$ is produced, i.e, the binary variable $z_p$ should be one whenever $C_{A,p}^{min} \leq C_A \leq C_{A,p}^{max}$ and zero in all other cases. 
As small variations between model \circled{1} and the nonlinear process model can occur, we add a small safety margin $\varepsilon_p$:
\begin{align}
    \tag{S1}
     -(1-z_p) \leq C_A - C_{A,p}^{min} - \varepsilon_p ~~~~~ \forall p \label{eq:prod_band_1} \\
     \tag{S2}
     -(1-z_p) \leq -C_A + C_{A,p}^{max} - \varepsilon_p~~~~~ \forall p \label{eq:prod_band_2}
\end{align}
We choose $\varepsilon_p = 0.003~\frac{mol}{L}$ based on simulation results. If the process cannot be simulated, a safety margin needs to be chosen based on measurement data. It would be reasonable to first start with a conservative safety margin, evaluate the difference between model and plant data and then adjust the safety margin.
Similar to the production binaries $z_p$, we need constraints that enforce the correct behavior of the binary $z_{cool}^{steady}$ needed for the cooling model (compare to equation (28) of the main paper):
\begin{align}
\tag{S3}
    & C_A - 0.29 \geq -(1-z_{cool}^{steady}) \label{eq:cool_steady_1}\\
    \tag{S4}
    & C_A - 0.31 \geq z_{cool}^{steady} \label{eq:cool_steady_2}
\end{align}
Lastly, we need to make sure that the active chillers have more spare capacity than the maximum error of the data-driven model.
As discussed in the Method Section of the main paper, our assumption that errors in energy demand are compensated by the energy system is only valid if spare capacity is given.
Consequently, we want to avoid a situation in which the scheduling optimization drives all active compression chillers to full load. 
If for example two chillers are active and operate at full load and the true nonlinear cooling power demand is slightly higher than the cooling power calculated from model \circled{2}, an unplanned start-up of chiller three has to be performed.
However, we want to avoid excessive on/off-switches as they shorten the equipment life-time. 
To this end, based on the fit of model \circled{2}, we introduce a safety margin of $\varepsilon_Q = 0.1$ by setting 
\begin{align}   
    \tag{S5}
    & Q_{cool} \leq (1-\varepsilon_Q)\sum_{i=1}^{3} z_{on,CC,i} Q_{CC,1}^{max}, \label{eq:safety_on}
\end{align}
such that there always is a 10~\% safety margin between the amount of cooling the active chillers can supply at full load and the cooling power calculated from model \circled{2}. Again, this safety margin could alternatively be determined based on measurements.

To speed up the solution time, we add the constraint that every product can only be produced once during the 24-hour horizon. 
Note that this additional constraint does not cut off the optimal solution: As all the transition times are between 0.5 and 1 hour, it is not reasonable to perform more transitions than necessary, because any additional transition would correspond to at least half an hour without product revenue, which clearly outweighs energy cost savings.
The constraint limiting the number of production starts can be seen as a discrete time analogue to setting the number of production slots equal to the number of products in a continuous time formulation, see, e.g., \cite{FloresTlacuahuac.2006}.
To implement the constraint, we introduce a continuous variable $s_{p,t}$ that is equal to one, whenever production of product $p$ starts in a timestep $t$, and zero otherwise.
Here, we use the index $t$ to indicate the timestep. 
For the sake of readability, this index is neglected in all the previous equations as the timestep is the same for all variables there. For every discrete timestep $t$, the following equations from \cite{VerenaNeisen.2018} are implemented such that $s_{p,t}$ is one if and only if the production binary $z_{p,t}$ is one and the production binary $z_{p,t-1}$ is zero:
\begin{align}
    \label{eq:start_1}
    \tag{S6}
    &0 \leq s_{p,t} \leq 1 \\
    \tag{S7}
    &s_{p,t} \geq z_{p,t} - z_{p,t-1} \\
    \tag{S8}
    &s_{p,t} \leq 1 - z_{p,t-1} \\
    \tag{S9}
    &s_{p,t} \leq z_{p,t} \\
    \tag{S10}
    \label{eq:limit_start}
    &\sum_{t \in \mathbb{T}_{dis}} s_{p,t} \leq 1 ~~~~~ \forall p
\end{align}
Note that one can easily change \cref{eq:limit_start} to allow two or more production starts if desired.

To formulate the optimization problem, we use pyomo \citep{hart2017pyomo,hart2011pyomo} and the extension pyomo.dae \citep{Nicholson.2018} to discretize time.
We solve the problem using gurobi\citep{gurobi} version 9.1.2  with an optimality gap of 1.0~\%. 
All calculations are performed on a Windows 10 machine with Intel(R) Core(TM) i5-8250U core and 24 GB RAM.

\subsection*{Details on the sequential steady-state scheduling benchmark}
In the sequential steady-state benchmark, we first optimize the production process decisions and afterwards the operation of the energy system, i.e., the compression chillers.

The optimization of the production process is a steady-state scheduling with predefined minimum transition times between the products using the pre-optimized transition profiles from the set-point filter tuning in the main paper.
For such an optimization, only the discrete time grid is necessary.
At each point in time, the process is either in steady-state production or in transition:
\begin{align}
    \label{eq:SEQ_status}
    \tag{S11}
    \sum_{p \in \mathbb{P}} z_p + \sum_{s \in \mathbb{S}} z_s = 1
\end{align}
In \cref{eq:SEQ_status}, $z_s$ is a binary indicating if the process is in transition $s$ with $\mathbb{S}$ being the set of possible transitions.

As the reactor cannot jump between products, a transition $s$ must be at least as long as the minimum transition time  $\Delta t_{s}^{min}$.
For implementation, we introduce the variable $\Delta t_s$ that indicates the time the process has been in a transition state $s$.
We show the relevant equations for transition $s=\RNum{1} \to \RNum{2}$. 
Again, the index $t$ is used to indicate the time step:
\begin{align}
\tag{S12}
    &\Delta t_{\RNum{1} \to \RNum{2},t} \leq z_{\RNum{1} \to \RNum{2},t} \\
    \tag{S13}
    &\Delta t_{\RNum{1} \to \RNum{2},t} \leq \Delta t_{\RNum{1} \to \RNum{2},t-1} + \Delta t_{dis}
\end{align}
The process can only be in the transition $s=\RNum{1} \to \RNum{2}$, if it already was in $s=\RNum{1} \to \RNum{2}$ in the last time step or if it was producing $\RNum{1}$:
\begin{align}
\tag{S14}
    z_{\RNum{1} \to \RNum{2},t} \leq z_{\RNum{1} \to \RNum{2},t-1} + z_{\RNum{1},t-1}
\end{align}
The process can only produce product $\RNum{2}$, if it was producing $\RNum{2}$ in last time step, or of it was in transition state $\RNum{1} \to \RNum{2}$ longer than the minimum time required for the transition ($\Delta t_{\RNum{1} \to \RNum{2},t-1} \geq \Delta t_{\RNum{1} \to \RNum{2}}^{min}$), or if it was in transition state $\RNum{3} \to \RNum{2}$ longer than the minimum time ($\Delta t_{\RNum{3} \to \RNum{2},t-1} \geq \Delta t_{\RNum{3} \to \RNum{2}}^{min}$). To implement this logic, we use big-M constraints \citep{Bemporad.1999}:
\begin{align}
    \label{eq:big_M_1}
    \tag{S15}
    z_{\RNum{2},t} \leq & 3z_{\RNum{2},t-1}M + (\Delta t_{\RNum{1} \to \RNum{2},t-1} - \Delta t_{\RNum{1} \to \RNum{2}}^{min})M +  3(1-z_{\RNum{2},t})M + 2z_{\RNum{3} \to \RNum{2},t-1}M \\
    \tag{S16}
    \label{eq:big_M_2}
    z_{\RNum{2},t} \leq& 3z_{\RNum{2},t-1}M + (\Delta t_{\RNum{3} \to \RNum{2},t-1} - \Delta t_{\RNum{3} \to \RNum{2}}^{min})M +   3(1-z_{\RNum{2},t})M + 2z_{\RNum{1} \to \RNum{2},t-1}M
\end{align}
In \cref{eq:big_M_1,eq:big_M_2}, $M$ needs to be a sufficiently large constant. We choose $M =2.4 \cdot 10^5$.
        
In the subsequent energy system optimization, the compression chillers are scheduled while the process schedule and thus the cooling demand are fixed.
Consequently, the optimization needs to decide on the on/off status and the load for each of the three chillers.
The chillers are modeled using equations (30) and (31) from the main paper and the minimum part-load constraint is implemented using equation (10) from the main paper.
Moreover, the energy balance equation (32) from the main paper is included, which sets the output of the three chillers equal to the cooling power $Q_{cool}$.
The cooling power $Q_{cool}$ is assumed to be equal to the steady-state cooling power $Q_{cool,p}^{steady}$ of a product $p$ during the production phases.
During transition phases, we use the average cooling power $Q_{cool,s}^{average}$ for the respective transition $s$. 
This average power is calculated from the simulations shown in Figure 7 of the main paper by dividing the integral of the cooling power by the transition length.
However, during a transition, the peak cooling demand is significantly larger than the average cooling demand.
Consequently, the active chillers have to provide enough spare capacity to be able to supply the peak cooling demand $Q_{cool,s}^{peak}$ during a transition $s$.
The calculation of the cooling power is given by \cref{eq:cool_bench} and the constraint to enforce spare capacity by \cref{eq:peak_sequential}:
\begin{align}
\label{eq:cool_bench}
    \tag{S17}
    Q_{cool} = \sum_{p \in \mathbb{P}} z_p Q_{cool,p}^{steady} + \sum_{s \in \mathbb{S}} z_s Q_{cool,s}^{average} \\
    \tag{S18}
    \label{eq:peak_sequential}
    \sum_{s \in \mathbb{S}} z_s Q_{cool,s}^{peak} \leq \sum_{i=1}^{3} z_{on,CC,i} Q_{CC,i}^{max}
\end{align}
Note that the binaries $z_s$ and $z_p$ are  fixed in the compression chiller scheduling as these decisions have already been taken by the production process scheduling.

Both optimization problems in the sequential benchmark, i.e., the production process scheduling and the subsequent compression chiller scheduling, are MILPs. 

\section*{Case study 3: Distillation column with variable purity}

\subsection*{Distillation column model}

Here, we briefly summarize the model equations for the distillation column taken from \cite{Skogestad.1988} together with a liquid flow model from \cite{Skogestad.1990}. 
The derivative of a liquid mole fraction $x_k$ on a tray $k$ that is not reboiler ($k=1$), feed ($k=N_F$), or condenser ($k=N+1$) is given by:
\begin{align}
    \tag{S19}    
    &\frac{dx_k}{dt} = \frac{1}{M_k} (L_{k+1} x_{k+1} + V y_{k-1} - L_k x_{k} - V y_{k}  )\\ \nonumber &~~~~~~~~~~~~~~\forall k =2,...,N_F-1 \text{ and } \forall k =N_F+1,...,N
\end{align}
where $M_k$ is the liquid hold-up, $L_k$ is the liquid flow leaving tray $k$ and $y_k$ the vapor mole fraction which is a function of the corresponding liquid mole fraction, i.e.,
\begin{align}
\tag{S20}
    y_{k} = \frac{\alpha x_{k}}{1 + (\alpha -1)x_{k}}~~\forall k =1,...,N,
\end{align}
with relative volatility $\alpha$.
Additionally, the total mass balance gives the derivative of the liquid hold-up, i.e.,
\begin{align}
    \tag{S21}
    \frac{M_k}{dt} = L_{k+1}  - L_k ~~\forall k =2,...,N_F-1 \text{ and } \forall k = N_F+1,...,N,
\end{align}
and the liquid flow $L_k$ is given as
\begin{align}
    \tag{S22}
    L_k = L_{k,0} + \frac{M_k - M_{k,0}}{\tau_k},
\end{align}
where $L_{k,0}$ and $M_{k,0}$ are liquid flow and hold-up at nominal conditions and $\tau_k$ is the hydraulic lag constant.
The equations for reboiler ($k=1$), feed ($k=N_F$), and condenser ($k=N+1$) are given by:
\begin{align}
    \label{eq:x1}
    \tag{S23}
    &\frac{dx_1}{dt} = \frac{1}{M_1} (L_2 x_2 - V y_1 - B x_1) \\
    \label{eq:M1}
    \tag{S24}
    &\frac{M_1}{dt} = L_2 - V - B \\
    \tag{S25}
    &\frac{dx_{N_F}}{dt} = \frac{1}{M_{N_F}} (L_{N_F+1} x_{N_F +1 } + V y_{N_F-1} + F z_F - L_{N_F} x_{N_F} - V y_{N_F}  ) \\
    \tag{S26}
    &\frac{M_{N_F}}{dt} = L_{N_F+1} + F  - L_{N_F} \\
    \tag{S27}
    &\frac{dx_{N+1}}{dt} = \frac{1}{M_{N+1}} (V y_{N}  - L x_{N+1} - D x_{N+1}) \\
    \tag{S28}
    &\frac{dM_{N+1}}{dt} = V - L - D
\end{align}
For further discussion of the model equations and the underlying assumptions, we refer to the original publications \citep{Skogestad.1988,Skogestad.1990}.
The parameter values are given in Table~\ref{tab:column_parameter}. The assumed limits on the inputs $u$ are given in Table~\ref{tab:input_bounds}.

\begin{table}[h]
    \centering
    \caption{Distillation column model parameters from the generic benchmark column A \citep{Skogestad.1988,Skogestad.1990}.}
    \large
    \renewcommand\theadfont{\large}
    \label{tab:column_parameter}
    \begin{tabular}{lc}
    \hline
    \thead{symbol} & \thead{value} \\
    \hline
    $z_F$ & $0.5$  \\
    $\alpha$ & 1.5 \\
    $N$ & 40 \\
    $N_F$ & 21 \\
    $M_{k0} / F$ & 0.5 min\\
    $\tau_k$ & 0.063 min \\
    \hline  
    \end{tabular}
\end{table}   

\begin{table}[h!]
    \centering
    \caption{Input bounds for distillation column model, all values are normalized to the feed flow F.}
    \large
    \renewcommand\theadfont{\large}
    \label{tab:input_bounds}
    \begin{tabular}{lrr}
     \hline
    input & lower bound & upper bound  \\
     \hline
    $D$ & 0  & 1 \\
    $B$ & 0 & 1\\
    $V$ & 1.5  & 2.5 \\
    $L$ & 1 & 2\\
     \hline
    \end{tabular}
\end{table}
    
\subsection*{Energy system component models}
For the two CHPs, the model from \cite{Voll.2013} is used, which gives the nominal thermal (th) and electric (elec) efficiencies $\eta$ as functions of the nominal heat output $Q^{nom,th}$:
\begin{align}
    \tag{S29}
    \eta^{nom,th} = -3.55\cdot 10^{-5} \frac{Q^{nom,th}}{1kW} + 0.498 \\
    \tag{S30}
    \eta^{nom,elec} = 3.55\cdot 10^{-5} \frac{Q^{nom,th}}{1kW} + 0.372 
\end{align}
The part load efficiencies are given as a functions of the heat output $Q^{th}$:
\begin{align}
\tag{S31}
    &\eta^{th}(Q^{th}) = \eta^{nom,th} \left(1.0960 - 0.0199q -0.0768q^2 \right) \\
    \tag{S32}
    &\eta^{elec}(Q^{th}) = \eta^{nom,elec} \left(0.5868 + 0.6743q -0.2611q^2 \right) \\ \nonumber
    &\text{with~~} q = \frac{Q^{th}}{Q^{nom,th}}
\end{align}
For the electrically-driven boiler, a constant thermal efficiency of 99~\% is assumed \citep{Baumgartner.2020}. 
As mentioned in the main paper, the minimum part-load is 50~\% for the CHPs and 20~\% for the electrically-driven boiler.

\subsection*{General coupling function}
In our case study, we couple the values of the top and bottom purities, $y_D$ and $x_B$, respectively.
This coupling needs to be chosen such that the ratio of top product flow rate $D$, and bottom flow rate $B$ stays constant in steady state.
Otherwise the flexible operation would change the ratio of top and bottom product.
In our case study, both top flow $D$ and bottom flow $B$ are equal to half of the feed flow $F$ in steady state.
To calculate the coupling between $y_D$ and $x_B$ that maintains this ratio, we consider the mass balances around the column in steady state:
\begin{align}
    \label{eq:mass_balance_comp}
    \tag{S33}
    z_F F &= y_D D + x_B B \\
    \tag{S34}
    \label{eq:mass_balance_tot}
    F &= D + B
\end{align}
Solving \cref{eq:mass_balance_tot} for $D$ and replacing $D$ in \cref{eq:mass_balance_comp} gives:
\begin{align}
    \label{eq:x_B_of_y_D}
    \tag{S35}
    &z_F F = y_D (F-B) + x_B B \\
    \nonumber
    \Leftrightarrow ~~& x_B = \frac{F}{B}z_F - y_D\left(\frac{F}{B}-1\right)
\end{align}
In our case study, we have $z_F=0.5$ (Table~\ref{tab:column_parameter}) and $\frac{F}{B}=2$ in steady state.
Thus, \cref{eq:x_B_of_y_D} simplifies to $x_B = 1- y_D$ (compare to equation~(36) in the main paper).

\subsection*{PI-controllers}
The PI controller equations are given by
\begin{align}
    \tag{S36}
    &V = K_{P,V} \left(e_{x_B} + \frac{1}{\tau_{I,V}} \int_0^t e_{x_B} dt\right) +V^{PI,0}~~~&\text{, with}~ e=(1-\rho_{SP,fil}) - x_B, \\
    \tag{S37}
    &L = K_{L,V} \left(e_{y_D} + \frac{1}{\tau_{I,L}} \int_0^t e_{y_D} dt\right) +L^{PI,0}~~~&\text{, with}~ e=\rho_{SP,fil} - y_D,
\end{align}
where $V^{PI,0}$ and $L^{PI,0}$ are steady-state values at nominal purity $\rho^{nom}$. 
The tuning parameters are tuned manually and set to $K_{P,V} = - 792/F$, $\tau_{I,L} = 25.8$~min, $K_{L,V} = 838/F$, and $\tau_{I,L} = 27.3$~min. 
The stable and accurate tracking is shown in Figure~\ref{fig:transitions_colum}.

\subsection*{Set-point filter tuning}

\begin{figure}   
    \centering
   \includegraphics{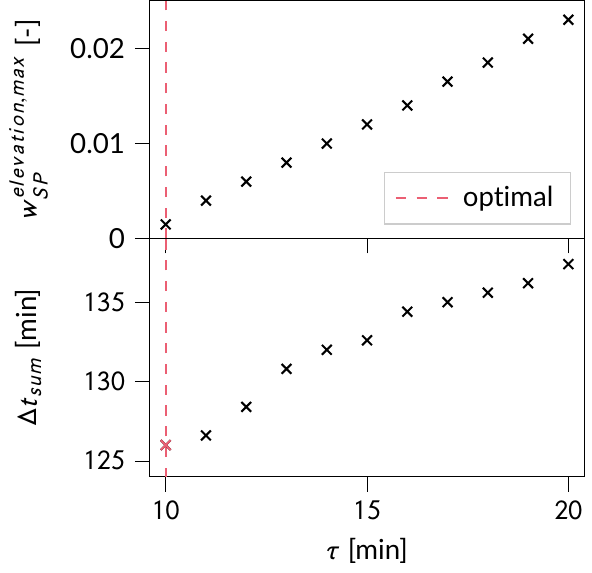}
    \caption{Result of the parameter tuning for the third case study: maximum allowable set-point elevation $w_{SP}^{elevation,max}$ (top) and sum of all 6 transition times $\Delta t_{sum}$ resulting with $w_{SP}^{elevation,max}$ (bottom) for different values of the time constant $\tau$. We choose the optimal, i.e., smallest $\Delta t_{sum}$, by setting $\tau=10$~min and $w_{SP}^{elevation}=0.0015$.}
   \label{fig:fit_tau}
\end{figure}

\begin{figure}
    \centering
    \includegraphics{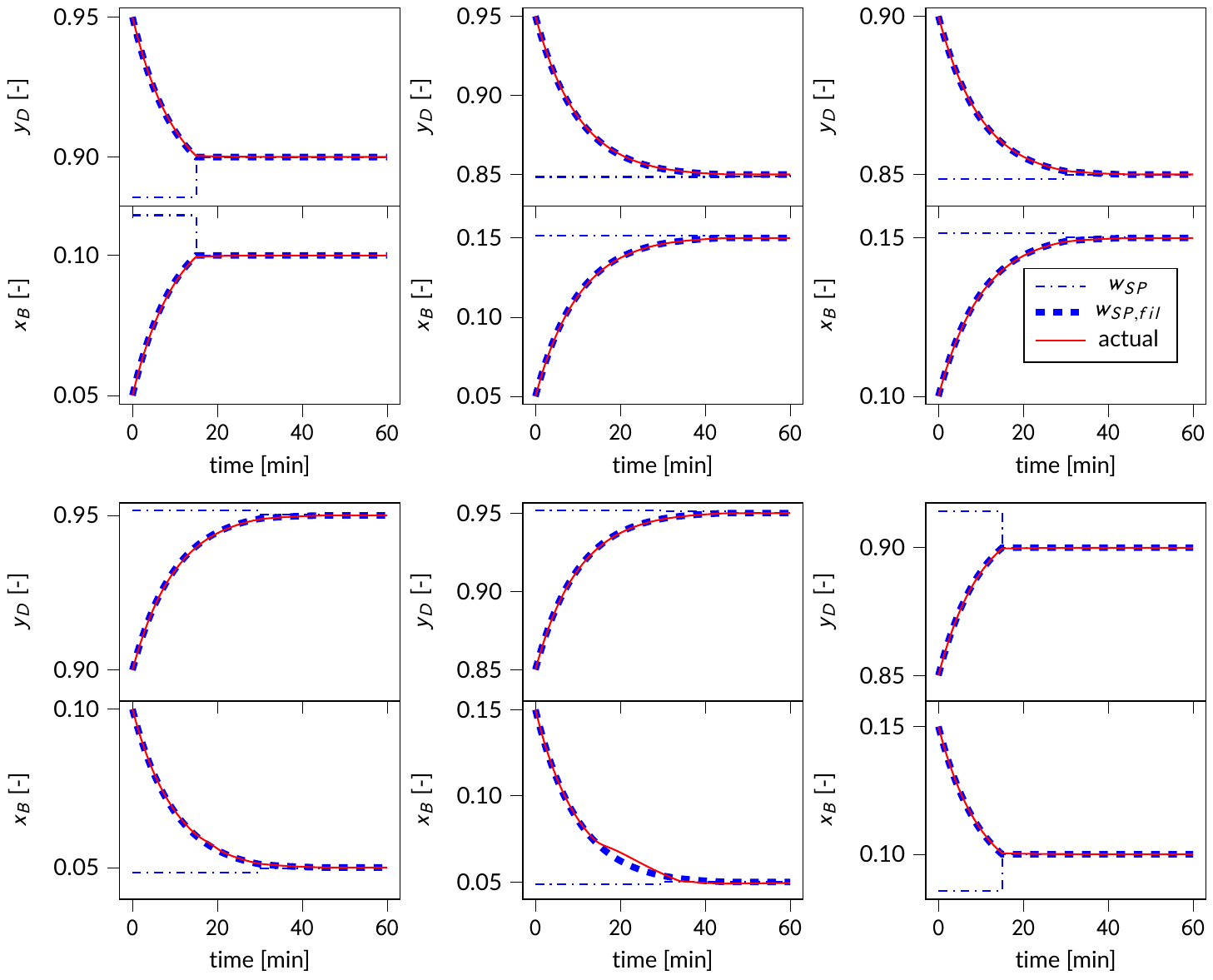}
    \caption{Six possible transitions between the three operating points $\rho^{max} = 0.95$, $\rho^{nom} = 0.9$, $\rho^{min} = 0.85$ in the third case study, with the chosen time constant $\tau=10~$min and set-point elevation  $w_{SP}^{elevation}=0.0015$ (compare to Figure S1 of the SI). The piece-wise constant set-point from optimization $w_{SP}$ results in a filtered set-point $w_{SP,fil}$, which can be tracked accurately such that the actual mole fractions $y_D$ and $x_B$ resulting from the controlled nonlinear process model match the filtered set-point $w_{SP,fil}$ well.}
   \label{fig:transitions_colum}
\end{figure}

To tune the parameter of the set-point filter, we consider the three operating points $\rho^{max} = 0.95$, $\rho^{nom} = 0.9$, $\rho^{min} = 0.85$ and the six possible transitions between these operating points.
In principle, we follow the same procedure as discussed in the first case study in the main paper.
However, as there are no products bands in this case study, we define a transition to be completed if both mole fractions $y_D$ and $x_B$ reach their final value within a tolerance of $\varepsilon = 0.005$.
Moreover, a transition is defined to be feasible if the maximum difference between a mole fraction and its filtered set-point $w_{SP,fil}$ is smaller than $\varepsilon = 0.005$.
Starting with $\tau = 1$~min and $w_{SP}^{elevation}=0$, we increase $\tau$  by $\Delta \tau = 1$~min and find that $\tau = 10$~min is the first time constant where all transitions are feasible.
We continue to increase $\tau$  by $\Delta \tau = 1$~min and for every $\tau$ we increase the allowed set-point elevation by steps of size $\Delta w_{SP}^{elevation} = 0.0005$ until one of the transitions becomes infeasible.
In contrast to case study  1, the smallest feasible time constant $\tau = 10$~min leads to the smallest sum of all six transition times $\Delta t_{sum}$  (Figure~\ref{fig:fit_tau}).
The 6 resulting transitions with the chosen values are shown in Figure~\ref{fig:transitions_colum} for both mole fractions $y_D$ and $x_B$.

\subsection*{Results}
\begin{figure}
    \centering
    \includegraphics[width=15cm]{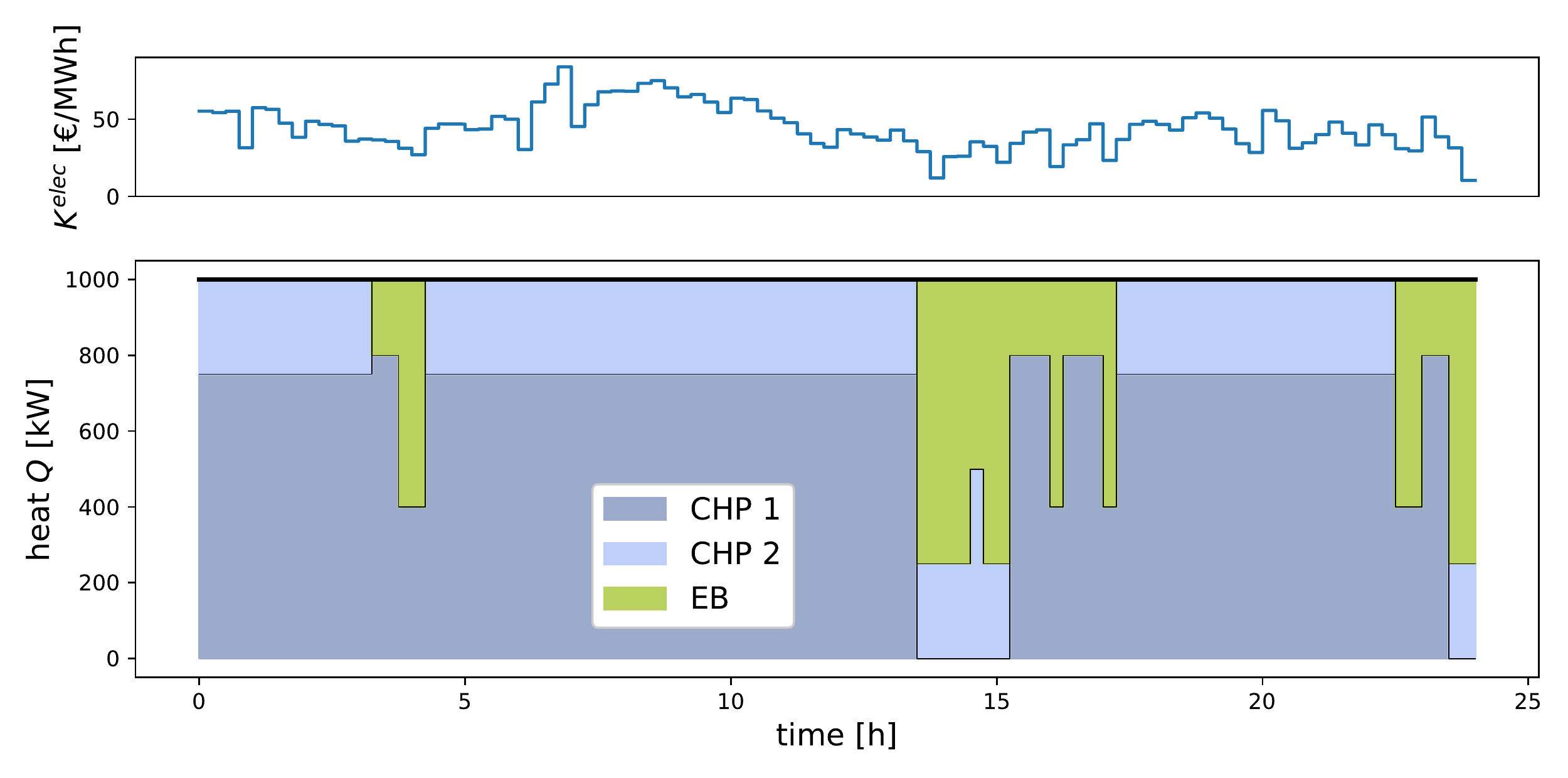}
    \caption{Operation in the steady-state benchmark for the third case study. Top: Electricity price $K^{elec}$. Bottom: The heat demand $Q$ of the column (bold black line) is shown and the portions of the heat demand supplied by the two combined heat and power plants (CHP1 and CHP2) and the electrically-driven boiler (EB) are indicated with colors. }
    \label{fig:case_study_3_steady_state_operation}
\end{figure}

\begin{figure}[h] 
    \centering
    \includegraphics{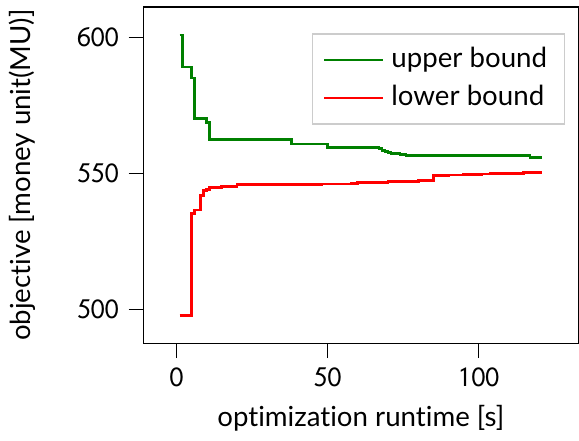}
    \caption{Convergence plot of scheduling optimization in the third case study.}
    \label{fig:convergence}
\end{figure}

In this Section, supplementary results for the third case study, are shown: The operation of the steady-state benchmark is shown in Figure~\ref{fig:case_study_3_steady_state_operation} (compare to Figure 16 of the main paper). 
The optimization favors the electrically-driven boiler for low power prices and operation of the CHPs for high power prices. 
Additionally, the convergence plot of our simultaneous dynamic scheduling approach is given in Figure~\ref{fig:convergence}. 
A near optimal solution with only 2.2~\% gap to the final lower bound is found after 11 seconds which is less than one tenth of the total runtime.

\twocolumn
\singlespacing
\section*{Nomenclature} \label{sec:nomenclature}
\noindent\textbf{Abbreviations} \\
\noindent
\begin{tabularx}{\columnwidth}{lX}
    CC & compression chiller \\
    CHP & combined heat and power plant \\
    CSTR & continuous-stirred-tank-reactor \\
    EB & electrically-driven boiler\\
    MILP & mixed-integer linear programming \\
    MU & money unit \\
    PI & proportional–integral
\end{tabularx} \\

\noindent\textbf{Greek symbols} \\
\noindent
\begin{tabularx}{\columnwidth}{lX}
    $\alpha$ & relative volatility \\
    $\varepsilon$ & safety margin \\
    $\eta$ &efficiency \\
    $\varrho$ & density\\
    $\rho$ & scheduling-relevant variable \\
    $\tau$ & time constant \\
\end{tabularx} \\

\noindent\textbf{Latin symbols} \\
\noindent
\begin{tabularx}{\columnwidth}{lX}
    $B$ & bottom flow rate \\
    $C_A$ & concentration of component A \\
    $C_P$& heat capacity \\
    $D$ & distillate flow rate \\
    $E_A$ & activation energy \\
    $F$ & feed flow rate \\
    $\Delta H_r$ & enthalpy of reaction \\
    $K$&price \\
    $K_P$ & proportional controller constant \\ 
    $k$ & reaction constant \\
    $L$ & reflux flow rate \\
    $L_k$ & liquid flow rate on tray $k$\\
    $M$ & big-M constant \\
    $M_k$ & hold-up of tray $k$ \\
    $N$ & number of trays \\
    $N_{cp}$ & order of collocation polynomial \\
    $Q$ & thermal power \\
    $q$ & flow rate \\
     $R$ &gas constant \\
    $s$ & variable to indicate start of production \\
    $T$ & temperature \\
    $t$ & time \\
    $V$ & volume (case study 1\&2) or boilup flow rate (case study 3) \\
    $w_{SP}$ & set-point \\
    $x_k$ & liquid mole fraction on tray k \\
    $y$ & continuous variable \\
    $y_k$ & vapor mole fraction on tray k \\
    $z$ & discrete variable \\
    $z_F$ &  mole fraction of feed flow 
\end{tabularx} \\

\noindent\textbf{Sets} \\
\noindent
\begin{tabularx}{\columnwidth}{lX}
    $\mathbb{P}$ & products \\
    $\mathbb{S}$ &set of possible transitions \\
    $\mathbb{T}_{dis}$ & timepoints on discrete grid
\end{tabularx} \\

\newpage
\noindent\textbf{Subscripts} \\
\noindent
\begin{tabularx}{\columnwidth}{lX}
    $B$ & reboiler \\
    $cont$ & continuous \\
    $cool$ & cooling \\
    $D$ & condenser \\
    $dis$ & discrete \\
    $elec$ &electricity \\
    $fil$ & filtered\\
    $on$ & on-off status \\
    $p$ & product \\
     $s$ & transition  \\
   $sum$ & summed value
\end{tabularx} \\

\noindent\textbf{Superscripts} \\
\noindent
\begin{tabularx}{\columnwidth}{lX}
    $elec$ &electric \\
    $max$ & maximum value\\
    $min$ & minimum value\\
    $nom$ & nominal \\
    $steady$ & steady-state value \\
    $th$ & thermal 
\end{tabularx} \\

\onecolumn
\bibliographystyle{apalike}
\renewcommand{\refname}{Bibliography} 
\bibliography{literature.bib}